\documentclass[leqno,11pt]{article}

\usepackage{amssymb,amsmath,amscd,amsfonts,a4,psfrag,graphicx,latexsym,epsfig,color}
\usepackage[all]{xy}

\usepackage{framed}

\usepackage{tikz}

\usepackage{verbatim}

\usetikzlibrary{matrix}

\newtheorem{theorem}{Theorem}[section]
\newtheorem{lemma}[theorem]{Lemma}
\newtheorem{proposition}[theorem]{Proposition}
\newtheorem{corollary}[theorem]{Corollary}
\newtheorem{definition}[theorem]{Definition}
\newtheorem{remark}[theorem]{Remark}
\newtheorem{example}[theorem]{Example}

\newtheorem{assumption}[theorem]{Assumption}

\newcommand{\agot}{\ensuremath{\mathfrak{a}}}

\newcommand{\ggot}{\ensuremath{\mathfrak{g}}}
\newcommand{\hgot}{\ensuremath{\mathfrak{h}}}
\newcommand{\kgot}{\ensuremath{\mathfrak{k}}}
\newcommand{\lgot}{\ensuremath{\mathfrak{l}}}
\newcommand{\mgot}{\ensuremath{\mathfrak{m}}}
\newcommand{\pgot}{\ensuremath{\mathfrak{p}}}

\newcommand{\ngot}{\ensuremath{\mathfrak{n}}}

\newcommand{\tgot}{\ensuremath{\mathfrak{t}}}
\newcommand{\ugot}{\ensuremath{\mathfrak{u}}}
\newcommand{\Rgot}{\ensuremath{\mathfrak{R}}}
\newcommand{\Sgot}{\ensuremath{\mathfrak{S}}}
\newcommand{\Xgot}{\ensuremath{\mathfrak{X}}}

\newcommand{\Ccal}{\ensuremath{\mathcal{C}}}

\newcommand{\Ecal}{\ensuremath{\mathcal{E}}}

\newcommand{\Fcal}{\ensuremath{\mathcal{F}}}

\newcommand{\Pcal}{\ensuremath{\mathcal{P}}}
\newcommand{\Rcal}{\ensuremath{\mathcal{R}}}

\newcommand{\Xcal}{\ensuremath{\mathcal{X}}}

\newcommand{\Zcal}{\ensuremath{\mathcal{Z}}}
\newcommand{\Ucal}{\ensuremath{\mathcal{U}}}
\newcommand{\Vcal}{\ensuremath{\mathcal{V}}}

\newcommand{\J}{\ensuremath{\mathbb{J}}}
\newcommand{\Z}{\ensuremath{\mathbb{Z}}}
\newcommand{\C}{\ensuremath{\mathbb{C}}}
\newcommand{\B}{\ensuremath{\mathbb{B}}}
\newcommand{\F}{\ensuremath{\mathbb{F}}}
\newcommand{\G}{\ensuremath{\mathbb{G}}}

\newcommand{\R}{\ensuremath{\mathbb{R}}}
\newcommand{\N}{\ensuremath{\mathbb{N}}}
\newcommand{\Pbb}{\ensuremath{\mathbb{P}}}

\newcommand{\croc}{\ensuremath{\hookrightarrow}}

\newcommand{\T}{\ensuremath{\hbox{\bf T}}}
\newcommand{\tr}{\ensuremath{\hbox{\bf Tr}}}

\newcommand{\Pol}{{\ensuremath{\rm Pol}}}

\def \tG {\widetilde{G}}
\def \tK {\widetilde{K}}
\def \tU {\widetilde{U}}

\def \tPbb {\widetilde{\Pbb}}
\def \tBbb {\widetilde{\B}}

\def \tggot {\tilde{\ggot}}

\def \tkgot {\tilde{\kgot}}
\def \tpgot {\tilde{\pgot}}
\def \tugot {\tilde{\ugot}}
\def \ttgot {\tilde{\tgot}}
\def \tngot {\tilde{\ngot}}
\def \tagot {\tilde{\agot}}

\def \ta {\tilde{a}}

\def \tk {\tilde{k}}
\def \tp {\tilde{p}}

\def \tg {\tilde{g}}

\def \tFcal{\widetilde{\mathcal{F}}}
\def \tX {\widetilde{X}}
\def \tY {\widetilde{Y}}

\def \tT {\widetilde{T}}

\def \txi {\tilde{\xi}}

\def \tw {\tilde{w}}

\newcommand{\horn}{\ensuremath{\hbox{\rm Horn}}}
\newcommand{\sing}{\ensuremath{\hbox{\rm Singular}}}

\begin{document}

\title{Moment polytopes in real symplectic geometry II : applications to singular value inequalities}

\author{Paul-Emile Paradan\footnote{IMAG, Univ Montpellier, CNRS, email : paul-emile.paradan@umontpellier.fr}}

\maketitle

\date{}


\begin{abstract}
In this work, we study some convex cones associated to isotropic representations of symmetric spaces. We explain the inequalities that describe them 
by means of cohomological conditions. In particular, we study the singular Horn cone which is the counterpart of the classical Horn cone, where the 
eigenvalues of Hermitian square matrices are replaced by the singular values of rectangular matrices.
\end{abstract}



\tableofcontents

\section{Introduction}

This paper is concerned with convexity properties associated to isotropy representations of symmetric spaces. Let $G/K$ be a Riemmanian symmetric space of the non-compact type and let $\pgot:=\T_e G/K$ be the isotropy representation of the compact Lie group $K$ (that we suppose connected). The $K$-orbits in $\pgot$ are parametrized by a closed cone $\agot_+$ contained in a maximal abelian subspace $\agot\subset \pgot$.

Suppose that $G/K\croc \tG/\tK$ is an embedding of Riemmanian symmetric spaces of the non-compact type. Thus $K$ is a closed subgroup of $\tK$ and we have an equivariant orthogonal projection $\pi :\tpgot\to \pgot$. The main purpose of this article is the description of the following convex cone
$$
\Pi_{\pgot}(\tK,K):=\left\{(\tilde{\xi},\xi)\in \tagot_+\times\agot_+ ; \ K\xi\,\subset \,\pi\big(\tK\tilde{\xi}\,\big)\right\}. 
$$

An efficient way to show that $\Pi_{\pgot}(\tK,K)$ is convex is to use the classical duality between the symmetric spaces of the non-compact type 
with those of the compact type. In our study, we suppose that $G\subset\tG$ are linear. In this context, there exists complex reductive groups 
$U_\C\subset \tU_\C$ equipped with an anti-linear involution $\sigma$ such that 
the connected component of the fixed point subgroups $(U_\C)^\sigma\subset (\tU_\C)^\sigma$ and $U^\sigma\subset \tU^\sigma$ are 
respectively $G\subset\tG$ and $K\subset \tK$ (see \S \ref{sec:kahler-hamiltonian}).

We consider maximal torus $T\subset U$ and $\tT\subset \tU$, that are invariant under $\sigma$, and such that the subspaces
$\tgot^{-\sigma}$ and $\ttgot^{-\sigma}$ are of maximal dimension. Then we can choose Weyl chambers $\tgot_+$ (resp. $\ttgot_+$) such that 
$\tgot^{-\sigma}\cap \tgot_+$ (resp. $\ttgot^{-\sigma}\cap \ttgot_+$) parametrizes the $K$-orbits in $\ugot^{-\sigma}$ 
(resp. the $\tK$-orbits in $\tugot^{-\sigma}$).

The convex cone
$$
\Pi(\tU,U):=\left\{(\tilde{\xi},\xi)\in \ttgot_+\times\tgot_+ ; \ U\xi\,\subset \,\pi\big(\tU\tilde{\xi}\,\big)\right\}. 
$$
has been the subject of numerous studies in recent decades. Let us mention the works of Horn \cite{Horn}, Klyachko \cite{Klyachko}, Belkale \cite{Belkale06} 
and Knutson-Tao \cite{Knutson-Tao-99}, when $U=U(n)$ and $\tU=(U(n))^s$. The case $\tU=(U)^s$ was studied by Belkale-Kumar \cite{BK1} 
and Kapovich-Leeb-Millson \cite{KLM-memoir-08}, and the general setting $U\subset \tU$ was considered by Berenstein-Sjamaar \cite{Berenstein-Sjamaar} 
and Ressayre \cite{Ressayre10,Ressayre-IF-11}.

\newpage

In particular, Ressayre obtained in \cite{Ressayre-IF-11} a complete description of 
the general\footnote{A face of $\Pi(\tU,U)$ is called {\em general} when it intersects the interior 
of the cone $\ttgot_+\times\tgot_+.$} faces of $\Pi(\tU,U)$. When $\tU=U^s$, he shows that the general faces of $\Pi(U^s,U)$ are parameterized by 
the set of $(\Pbb,\Xgot^\Pbb_{w_1},\ldots, \Xgot^\Pbb_{w_{s+1}})$ where $\Pbb$ is a standard parabolic subgroup of $U_\C$ and the $\Xgot^\Pbb_{w_k}$'s 
are  Schubert varieties of $U_\C/\Pbb$ such that $[\Xgot^\Pbb_{w_1}]\odot_0\cdots\odot_0[\Xgot^\Pbb_{w_{s+1}}]=[pt]$. Here $\odot_0$ is the 
Belkale-Kumar's product on the cohomology groups $H^*(U_\C/\Pbb,\Z)$.

The link between $\Pi(\tU,U)$ and $\Pi_{\pgot}(\tK,K)$ is given by a theorem of O'Shea and Sjamaar \cite{OSS99} that we recall in the next Section: it says that we have 
a natural identification
$$
\Pi_{\pgot}(\tK,K)\simeq \Pi(\tU,U)\cap \ttgot^{-\sigma}\times \tgot^{-\sigma}.
$$
From the above isomorphism, the convex cone $\Pi_\pgot(\tK,K)$ is completely determined since the equations of the cone $\Pi(\tU,U)$ are known. 
But by doing so, the inequalities describing $\Pi_\pgot(\tK,K)$ are in general highly overdetermined.

Our main contribution is a description of a smaller list of inequalities describing $\Pi_\pgot(\tK,K)$. Let us give an overview of our result when $\tG=G^s$ and $\tK=K^s$ with $s\geq 2$. 
A more precise version is explained in Sections \ref{section:main-result} and \ref{sec:levi}.

\begin{theorem}\label{theo:intro}The convex cone $\Pi_\pgot(K^s,K)$ can be described by a system of inequalities parameterized by the set of 
$(\Pbb,\Xgot^\Pbb_{w_1},\ldots, \Xgot^\Pbb_{w_{s+1}})$ where 
\begin{itemize}
\item $\Pbb$ is a $\sigma$-invariant parabolic subgroup of $U_\C$,
\item  $\Pbb$ is maximal among the $\sigma$-invariant parabolic subgroups of $U_\C$,
\item the $\Xgot^\Pbb_{w_k}=\overline{\B[w_k]}$ are Schubert varieties such that 
$\B[w_k]\cap(U_\C/\Pbb)^\sigma\neq\emptyset$,
\item the relation $[\Xgot^\Pbb_{w_1}]\odot_0\cdots\odot_0[\Xgot^\Pbb_{w_{s+1}}]=[pt]$ holds  in $H^*(U_\C/\Pbb,\Z)$.
\end{itemize}
\end{theorem}

A description of $\Pi_\pgot(K^s,K)$ was also obtained by Kapovich-Leeb-Millson in \cite{KLM-memoir-08}. 
We explain in \S \ref{sec:KLM} why our result is more accurate than theirs.

When working with the reductive group $G=U(p,q)$ and its maximal compact subgroup $K=U(p)\times U(q)$, the corresponding convex cone $\Pi_\pgot(K^2,K)$ 
is the singular Horn cone which we denote by $\sing(p,q)$. In \S \ref{sec:singular-horn}, we describe $\sing(p,q)$ by means of the inequalities of Theorem 
\ref{theo:intro}.

\subsection*{Acknowledgements}  
I would like to thank Mich\`ele Vergne and Velleda Baldoni for their precious help in the calculation of some Littlewood-Richardson coefficients.

\subsection{K\"{a}hler-Hamiltonian manifold with involution}\label{sec:kahler-hamiltonian}

The purpose of this section is to explain the theorem of O'Shea and Sjamaar \cite{OSS99} that we mentioned above.

\subsubsection*{Duality}
Let $\tG\subset {\rm GL}(n,\R)$ be a connected linear real reductive subgroup with maximal compact subgroup $\tK\subset {\rm SO}(n,\R)$. 
 Let $\tpgot$ be the intersection of the Lie algebra $\tggot$ of $\tG$ with the subspace of $n\times n$ symmetric matrices. Let $\tU\subset GL(n,\C)
 \subset GL(2n,\R)$ be the compact connected Lie group with Lie algebra  $\tugot=\tkgot\oplus i\tpgot$. Let us consider the involution $\sigma(A)=J A J^{-1}$ of  $GL(2n,\R)$ defined by $J=\begin{pmatrix}
0& -I_n\\I_n & 0\end{pmatrix}$.  We see that $\sigma$ commutes with the Cartan involution, leaves stable the subgroups $\tU$, and that the morphism $j:{\rm GL}(n,\R)\to {\rm GL}(2n,\R)$, $j(g)={\rm Diag}(g,g)$ defines an isomorphism between $\tK$ and the connected component of $\tU^\sigma$. 

The compact Lie groups $\tU\subset U(n)\subset SO(2n,\R)$ admit a complexification $\tU_\C\subset GL(2n,\C)$ (see \S III.8 in
\cite{Brocker-TomDieck}). We extend the involution $\sigma$ on $GL(2n,\R)$ to the conjugate linear involution $\sigma(g)=J \overline{g} J^{-1}$ on $GL(2n,\C)$. We see then that the morphism $j:{\rm GL}(n,\R)\to {\rm GL}(2n,\C)$ defines an isomorphism between $\tG$ and the connected component of the subgroup $\tU_\C^\sigma$ fixed by the involution.

Let $\ttgot\subset \tugot$ be the Lie algebra of a maximal torus $\tT\subset \tU$ stable under the involution and such that 
$\tagot:=\frac{1}{i}\ttgot^{-\sigma}\subset\tpgot$ is of maximal dimension. One can choose the Weyl chamber $\ttgot_+\subset\ttgot$ appropriately so that the restricted Weyl chamber 
$\tagot_+=\frac{1}{i}\ttgot_+\cap \ttgot^{-\sigma}$ parameters the $\tK$-orbits on $\tpgot=\frac{1}{i}\tugot^{-\sigma}$ (see the Appendix in 
\cite{OSS99}).

If we work with a connected linear real reductive subgroup $G\subset \tG$ the same construction holds. We have a compact connected subgroup 
$U\subset \tU$ with Lie algebra  $\ugot=\kgot\oplus i\pgot$. Its complexification $U_\C\subset \tU_C$ is stable under the conjugate linear involution 
$\sigma$ and the real reductive group $G$ admits a natural identification with the connected component of the subgroup $U_\C^\sigma$. 
The maximal torus $T\subset U$ is taken invariant under $\sigma$, contained in $\tT$, and such that 
$\agot=\frac{1}{i}\tgot^{-\sigma}\subset \tagot$ is of maximal dimension. The Lie algebras $\ugot$ and $\tugot$ belong to the vector space 
$\mathfrak{gl}(2n,\R)$ that is equipped with the canonical Euclidean norm $\|X\|^2=\sqrt{{\rm Tr}(XX^t)}$. We denote by $\pi: \tugot \to \ugot$ the orthogonal projection.
 
\subsubsection*{O'Shea-Sjamaar's Theorem}

We consider the following geometrical setting : the group $\tU_\C$ viewed as a complex manifold equipped with the following action of 
$\tU_\C\times U_\C$ : $(\tilde{g},g)\cdot x= \tilde{g} x g^{-1}$. 

\medskip

We identify the complex reductive group $\tU_\C$ with the tangent bundle $\T\tU$ through the isomorphism 
$\varphi_2:\T\tU\to \tU_\C, \varphi(\ta,\tX)=\ta e^{i\tX}$ given by the Cartan decomposition. If we use $\varphi_2$ to transport 
the complex structure of $\tU_\C$ we get a complex structure $\J$ on $\T\tU$.

Since $\T^*\tU\simeq\tU\times \tugot^*$ and $\T\tU\simeq\tU\times \tugot$ through left translations, the isomorphism 
$\varphi_1:\T\tU\simeq\T^*\tU$ is obtain by means of the identification $\tugot\simeq \tugot^*$ given by the invariant scalar product. 
If we use $\varphi_1$ to transport the canonical symplectic form of $\T^*\tU$, the resulting symplectic structure $\Omega$ on $\T \tU$
is invariant relatively to the $\tU\times U$-action, and is compatible with the complex structure $\J$ (see \cite{Hall97}, \S 3). 
The moment map relative to the 
$\tU\times U$-action on $(\T\tU,\Omega)$ is a proper map $\Phi=\Phi_{\tU} \oplus \Phi_U : \T\tU \to \tugot\times \ugot$ defined by
\begin{equation}\label{eq:momentcotangent}
\Phi_{\tU}(\ta,\tX)=-\ta\tX,\quad \Phi_U(\ta,\tX)=\pi(\tX), 
\end{equation}

The conjugate linear involution $\sigma$ on $\tU_\C$ defines (through $\varphi_2$) an involution $\tau:\T\tU\to \T\tU$ defined by $\tau(\ta,\tX)=(\sigma(\ta),-\sigma(\tX))$.

The $\tU_\C\times U_\C$-manifold $\T\tU$ equipped with the data $(\Omega,\J,\Phi,\tau,\sigma)$ is a K\"ahler-Hamiltonian manifold with involution : 
\begin{itemize}
\item $(\T\tU,\Omega,\J)$ is a K\"ahler manifold, 
\item the $\tU_\C\times U_\C$-action preserves the complex structure $\J$,
\item the $\tU\times U$-action is Hamiltonian with proper moment map $\Phi$,
\item $\tau^*(\Omega)=-\Omega$ and $\tau^*(\J)=-\J$,
\item $\tau(g\cdot x)=\sigma(g)\cdot\tau(x)$, and $\Phi(\tau(x))=-\sigma(\Phi(x))$, $\forall(g,x)\in \tU_C\times U_C\times \T\tU$. 
\end{itemize}

We finish the section by recalling the  O'Shea-Sjamaar's Theorem \cite{OSS99}. Let $\Delta_{\tU\times U}(\T\tU)\subset \tilde{\tgot}_{+}\times \tgot_+$ be the Kirwan polytope associated to $\Phi$. Equations (\ref{eq:momentcotangent}) show that
$$
\Delta_{\tU\times U}(\T\tU)=\left\{(\tX,X)\in\tilde{\tgot}_{+}\times \tgot_{+}\, |\, -UX\subset \pi\big(\tU\tX\,\big) \right\}.
$$

We consider now the submanifold $(\T\tU)^\tau$ fixed by $\tau$. Let $\Zcal$ be the connected component of $(\T\tU)^\tau$ containing $(e,0)$: 
through the diffeomorphism 
$\varphi_2:\T\tU\simeq \tU_\C$, we see that $\Zcal\simeq \tG$. The restriction of the moment $\Phi$ to $\Zcal$ takes value in $\tugot^{-\sigma}\times \ugot^{-\sigma}$. If we use 
the canonical isomorphism $(\tX,X)\in \tpgot\times\pgot\longrightarrow i\,(\tX,X)\in \tugot^{-\sigma}\times \ugot^{-\sigma}$, we obtain a gradient map 
\begin{eqnarray*}
\Psi : \tG\simeq \tK\times\tpgot&\longrightarrow & \tpgot\times\pgot\\
(\tk,\tX)&\longmapsto & (-\tk\tX,\pi(\tX)).
\end{eqnarray*}
The corresponding Kirwan polytope is $\Delta_{\tK\times K}(\tG)=\{(\tilde{\xi},\xi)\in \tagot_+\times\agot_+ ; \ -K\xi\,\subset \,\pi\big(\tK\tilde{\xi}\,\big)\}$.

\begin{theorem}[O'Shea-Sjamaar]
Through the isomorphism $(\tX,X)\mapsto i(\tX,X), \tagot_+\times\agot_+\to \tilde{\tgot}_{+}\times \tgot_{+}\bigcap \ttgot^{-\sigma}\times \tgot^{-\sigma}$ we have 
an identification
$$
\Delta_{\tK\times K}(\tG)\simeq\Delta_{\tU\times U}(\T\tU)\bigcap \ttgot^{-\sigma}\times \tgot^{-\sigma}.
$$
\end{theorem}

If $w_0$ denotes the longuest element of the Weyl group $W=N_U(T)/T$, we see that the cones $\Delta_{\tU\times U}(\T\tU)$ and $\Pi(\tU,U)$ are related as follows: $(\tX,X)\in \Delta_{\tU\times U}(\T\tU)\Longleftrightarrow (\tX,-w_0 X)\in \Pi(\tU,U)$. The element $w_0$ is $\sigma$-invariant and satisfies $w_0(\agot_+)=-\agot_+$ (see \cite{OSS99}, Lemma 5.4). Hence we obtain the same relations :  $(\tX,X)\in \Delta_{\tK\times K}(\tG)\Longleftrightarrow (\tX,-w_0 X)\in \Pi_\pgot(\tK,K)$.

Finally, the O'Shea-Sjamaar's Theorem tells us that $\Pi_\pgot(\tK,K)\simeq \Pi(\tU,U)\bigcap \ttgot^{-\sigma}\times \tgot^{-\sigma}$. Nevertheless, 
this identity does not give an effective description of $\Pi_\pgot(\tK,K)$ because most of the inequalities of the cone  $\Pi(\tU,U)$ become redundant when we restrict them 
to $\ttgot^{-\sigma}\times \tgot^{-\sigma}$. In the next section, we propose a smaller list of inequalities describing $\Pi_\pgot(\tK,K)$.

\subsection{Description of the main result}\label{section:main-result}

Let $G\subset \tG\subset {\rm GL}(n,\R)$ be two connected linear real reductive subgroups with maximal compact subgroups $K\subset\tK\subset {\rm SO}(n,\R)$. Let $\agot\subset\tagot$ be maximal abelian subspaces of $\pgot\subset\tpgot$.

In this article, we work under the following assumption:

\begin{assumption}\label{hypothese}
If $\mathfrak{b}$ is an ideal of $\tggot$ contained in $\ggot$, then $\mathfrak{b}\subset Z(\tggot)\cap\ggot$.
\end{assumption}

\subsubsection{Admissible elements and polarized trace}

An element $\gamma\in \agot$ is called {\em rational} if the eigenvalues of the symmetric endomorphism ${\rm ad}(\gamma):\tggot\to\tggot$ are rational. 
We let $\Sigma(\tggot/\ggot)\subset \agot^*$ denote the set of non-zero weights relative to the 
$\agot$-action on $\tggot/\ggot$. 

\begin{lemma}
Under Assumption \ref{hypothese}, we see that ${\rm Vect}\left(\Sigma(\tggot/\ggot)\right)^\perp= Z(\tggot)\cap\pgot$.
\end{lemma}
{\em Proof :} The set $\mathfrak{b}:=\{X\in\ggot, [X,\tggot]\subset \ggot\}$ is an ideal of $\tggot$ contained in $\ggot$. Thanks to Assumption \ref{hypothese}, we see that 
$\mathfrak{b}=Z(\tggot)\cap\ggot$. Then the vector space ${\rm Vect}\left(\Sigma(\tggot/\ggot)\right)^\perp=\{X\in\agot, [X,\tggot]\subset \ggot\}$ coincides with 
$Z(\tggot)\cap\pgot=Z(\tggot)\cap\agot$. $\Box$

\medskip

If $\gamma\in\agot$, we denote by $\Sigma(\tggot/\ggot)\cap \gamma^\perp$ the subset of weights vanishing against $\gamma$.

\begin{definition}\label{def:admissible-exemple}
An element $\gamma\in \agot$ is {\em admissible} if it is {\em rational} and if 
\begin{equation}\label{eq:condition-gamma}
{\rm Vect}\big(\Sigma(\tggot/\ggot)\cap \gamma^\perp\big)={\rm Vect}\big(\Sigma(\tggot/\ggot)\big)\cap \gamma^\perp.
\end{equation}
\end{definition}

\medskip

If $\ell:E\to E$ is a symmetric endomorphism of an Euclidean space, we define 
$\tr_{\ell}(E^{\ell>0})=\sum_{a>0}a\,\dim(E^{\ell=a})$, where $E^{\ell=a}=\{v\in E, \ell(v)=av\}$.

If $\gamma\in\agot$, its adjoint action ${\rm ad}(\gamma):\ggot\to\ggot$ defines a symmetric endomorphism, hence for any subspace 
$\mgot\subset\ggot$ stable under ${\rm ad}(\gamma)$ one can define the number $\tr_{{\rm ad}(\gamma)}(\mgot^{{\rm ad}(\gamma)>0})$ that we denote simply by $\tr_{\gamma}(\mgot^{\gamma>0})$.

\subsubsection{Parabolic subgroups}

The maximal torus $T\subset U$ is stable under the involution $\sigma$, so the Weyl group $W=N_U(T)/T$ is equipped also with an involution still 
denoted $\sigma$, and we denote by $W^\sigma$ the subgroup fixed by it. Let $\Rgot$ be the set of roots relative to the action of $T$ on $\ugot_\C$.

The maximal compact subgroup $K\subset G$ correspond to the connected component of the subgroup $U^\sigma$. We consider the restricted Weyl group $W_\agot=N_K(\agot)/Z_K(\agot)$ relative to the maximal abelian subspace $\agot=\frac{1}{i}\tgot^{-\sigma}\subset\frac{1}{i}\ugot^{-\sigma}=\pgot$. Let $\Sigma$ be the set of (restricted) roots relative to the action of $\agot$ on $\ggot$.

Let us consider $\gamma\in\agot$ and the subgroups $W^\gamma\subset W$ and $W_\agot^\gamma\subset W_\agot$ fixing $\gamma$.
The following fact is explained in \S \ref{sec:bruhat}.
\begin{lemma}\label{lem:identification-weyl-group}
We have canonical identifications
$(W/W^\gamma)^\sigma\simeq W^{\sigma}/W^\sigma\cap W^\gamma\simeq W_\agot/W_\agot^\gamma$.
\end{lemma} 

We choose a Weyl chamber $\tgot_+\subset\tgot$ so that the restricted Weyl chamber $\agot_+=\frac{1}{i}\tgot_+\cap \tgot^{-\sigma}$ parameters the $K$-orbits on $\tpgot$. The choice of $\tgot_+$ determines a system of positive roots $\Rgot^+$, and  we denote by $\B\subset U_\C$ the Borel 
subgroup with Lie algebra $\tgot_\C\oplus \sum_{\alpha \in \Rgot^+}(\ugot_\C)_\alpha$.

To any element $\gamma\in \agot$, we associate the parabolic subgroup 
\begin{equation}\label{eq:parabolic-gamma}
\Pbb_\gamma=\{g\in U_\C, \lim_{t\to\infty}e^{t\gamma}ge^{-t\gamma}\ {\rm exists}\}.
\end{equation}
Notice that $\Pbb_\gamma$ is invariant under the involution $\sigma$ and that $\B\subset \Pbb_\gamma$ when $\gamma$ is anti-dominant, i.e. $\gamma\in -\agot_+$.

\subsubsection{Schubert Calculus}

Here we work with two complex reductive groups $U_\C\subset\tU_\C$ equipped with an antilinear involution $\sigma$. 
We choose maximal tori $T\subset\tT$ invariant by $\sigma$, such that  the corresponding subspaces $\agot:=\frac{1}{i}\tgot^{-\sigma}
\subset\tagot:=\frac{1}{i}\ttgot^{-\sigma}$ are maximal abelian in  $\pgot\subset\tpgot$.

To any $\gamma\in\agot$, we associate via (\ref{eq:parabolic-gamma}) the parabolic subgroups $\Pbb_\gamma\subset \tPbb_\gamma$ and 
we define the flag varieties  $\Fcal_\gamma= U_\C/\Pbb_\gamma$ and $\tFcal_\gamma= \tU_\C/\tPbb_\gamma$. Notice that the involution $\sigma$
defines an antiholomorphic involution on the flag varieties $\Fcal_\gamma$ and $\tFcal_\gamma$ that we still denoted by $\sigma$. As
$\Pbb_\gamma=U_\C\cap \tPbb_\gamma$, we have a canonical embedding $\iota:\Fcal_\gamma\croc \tFcal_\gamma$. We denote 
$\iota^*: H^{*}(\tFcal_\gamma,\Z)\to H^{*}(\Fcal_\gamma,\Z)$ the pullback in cohomology.

Thanks to the Bruhat decomposition, $\Fcal_\gamma=\bigcup_{w\in W/W^\gamma}\B[w]$, we know that the $\B$-orbits on $\Fcal_\gamma$ are parametrized by $W/W^\gamma$. We associate to any $w\in W/W^\gamma$, the Schubert cell $\Xgot^{o}_{w,\gamma}=\B [w]$, the Schubert variety
$\Xgot_{w,\gamma}:=\overline{\Xgot^{o}_{w,\gamma}}$ and its cycle class in cohomology $[\Xgot_{w,\gamma}]\in H^{*}(\Fcal_\gamma,\Z)$. 

Thanks to Lemma \ref{lem:identification-weyl-group}, we know that we can attach a Schubert cell $\Xgot^{o}_{w,\gamma}$ to any element 
$w\in W_\agot/W_\agot^\gamma\simeq (W/W^\gamma)^\sigma$. We will see in  \S \ref{sec:bruhat} that an element $w\in W/W^\gamma$ is $\sigma$-invariant  
if and only if $\Xgot^{o}_{w,\gamma}\cap(\Fcal_\gamma)^\sigma\neq\emptyset$.

\subsubsection{Main result}

Let $G\subset \tG$ be two linear real reductive groups satisfying Assumption \ref{hypothese}. Let $\pi:\tpgot\to\pgot$ be the orthogonal projection. 

Let $w_0\in W_\agot$ be the unique element such that $w_0(\agot_+)=-\agot_+$. The choice of Weyl chamber $\agot_+$ defines a system of positive roots $\Sigma^+$, and we denote by $\ngot=\sum_{\beta\in\Sigma^+}\ggot_\beta$ the corresponding real nilpotent subalgebra of $\ggot$.

Recall that the cohomology class $[pt]$ associated to a singleton $\{pt\}\subset \Fcal_\gamma$  is a basis of $H^{max}(\Fcal_\gamma,\Z)$.

\begin{theorem}\label{theo:main}  Let $(\tilde{\xi},\xi)\in\tagot_+\times \agot_+$. 
We have  $K\xi\subset \pi\big(\tK\txi\big)$ if and only if 
$$
(\tilde{\xi},\tilde{w}\gamma)\geq (\xi,w_0w\gamma)
$$
for any rational antidominant element $\gamma\in\agot$ and any $(w,\tilde{w})\in W_\agot/W_\agot^\gamma\times 
W_{\tagot}/W_{\tagot}^\gamma$ satisfying the following properties:
\begin{enumerate}
\item[a)] ${\rm Vect}\big(\Sigma(\tggot/\ggot)\cap \gamma^\perp\big)={\rm Vect}\big(\Sigma(\tggot/\ggot)\big)\cap \gamma^\perp$.
\item[b)] $[\Xgot_{w,\gamma}]\cdot \iota^*([\tilde{\Xgot}_{\tilde{w},\gamma}])= [pt]$ in $H^*(\Fcal_\gamma,\Z)$.
\item[c)] ${\rm Tr}_{w\gamma}(\ngot^{w\gamma>0})+{\rm Tr}_{\tw\gamma}(\tngot^{\tw\gamma>0})={\rm Tr}_{\gamma}(\tggot^{\gamma>0})$.
\end{enumerate}
The result still holds if we replace condition $b)$ by the weaker condition 
$$
b')\qquad [\Xgot_{w,\gamma}]\cdot \iota^*([\tilde{\Xgot}_{\tilde{w},\gamma}])=\ell [pt]\ {\rm with}\ \ell\geq 1,\  {\rm in}\  H^*(\Fcal_\gamma,\Z).
$$
\end{theorem}

\subsection{Bruhat decompositions}\label{sec:bruhat}

The complex parabolic subgroup $\Pbb_\gamma$ is stable under the involution $\sigma$. The intersection $P_\gamma:= \Pbb_\gamma\cap G\subset (\Pbb_\gamma)^\sigma$ defines a real parabolic subgroup of $G$ (that is not necessarily connected). Let us explain why the fixed point submanifold $(\Fcal_\gamma)^\sigma$ corresponds to the real flag variety $\Fcal_\gamma^\R=G/P_\gamma$. As the flag variety  $\Fcal_\gamma= U_\C/\Pbb_\gamma\simeq U/U^\gamma$, equipped with the involution $\sigma$, admits an identification with the 
adjoint orbit $U\gamma$, we have  
$$
(\Fcal_\gamma)^\sigma\simeq U\gamma \cap i\ugot^{-\sigma}= K\gamma\simeq K/K^\gamma\simeq G/P_\gamma=\Fcal_\gamma^\R.
$$ 
The crucial point here is the equality $U\gamma \cap i\ugot^{-\sigma}= K\gamma$ (see \cite{OSS99}, Example 2.9).

We fix an element $\gamma_o$ in the interior of the Weyl chamber $\agot_+$ and we consider the parabolic subgroup 
$\Pbb:= \Pbb_{-\gamma_o}\subset U_\C$. Then $P:=\Pbb\cap G$ is the minimal parabolic subgroup of $G$ with Lie algebra 
$Z_{\kgot}(\agot)\oplus \agot\oplus \ngot$.

Let us consider the following Bruhat decompositions (see \cite{DKV-83,Richardson92}) :
\begin{enumerate}
\item $\Fcal_\gamma=\bigcup_{w\in W/W^\gamma}\B[w]$ relative to the $\B$-action on $\Fcal_\gamma$,
\item $\Fcal_\gamma=\bigcup_{u\in W^{\gamma_o}\setminus W/W^\gamma}\Pbb[u]$ relative to the $\Pbb$-action on $\Fcal_\gamma$,
\item $\Fcal^\R_\gamma=\bigcup_{v\in W_{\agot}/W_{\agot}^\gamma} P[v]$ relative to the $P$-action on $\Fcal_\gamma^\R$.
\end{enumerate}

Now we give a proof of Lemma  \ref{lem:identification-weyl-group}. First we recall a standard result relating the Weyl groups $W_{\agot}$ and $W$ with the involution $\sigma$ (see the Appendix B in  \cite{OSS99}). The subgroup $W^\sigma$ fixed by $\sigma$ is equal to the normalizer subgroup $N_{W}(\agot)$. Since $\gamma_o$ is a regular element of $\agot$, the subgroup $W^{\gamma_o}$ coincides with the centralizer subgroup $Z_{W}(\agot)$, and the restricted Weyl group $W_{\agot}$ admits a canonical identification with $N_{W}(\agot)/Z_{W}(\agot)$.

Thus, we have two exact sequences 
$$
0\to W^{\gamma_o}\to W^\sigma\to W_{\agot}\to 0\qquad {\rm and}\qquad 0\to W^{\gamma_o}\to W^\sigma\cap 
W^{\gamma}\to W_{\agot}^\gamma\to 0,
$$
that induce an isomorphism $W^\sigma/W^\sigma\cap W^{\gamma}\simeq W_{\agot}/W_{\agot}^\gamma$.  The other isomorphism 
$(W/W^{\gamma})^\sigma\simeq W_\agot/W_\agot^\gamma$ follows from the fact that $W/W^{\gamma}\simeq W\gamma=U\gamma\cap i\tgot$ : it gives  
$$
(W/W^{\gamma})^\sigma\simeq (W\gamma)^\gamma=U\gamma\cap i\tgot^{-\sigma}=K\gamma\cap i\tgot^{-\sigma}= W_\agot\gamma
\simeq W_\agot/W_\agot^\gamma.
$$

Notice also that we have an inclusion $W^\sigma/W^\sigma\cap W^{\gamma}\croc W^{\gamma_o} \backslash W/W^\gamma$ : 
it is due to the fact that $W^{\gamma_o}$ is a distinguished subgroup of $W^\sigma$ contained in $W^{\gamma}$.
Finally we have proven $W_{\agot}/W_{\agot}^\gamma$ can be seen as a subset of 
$W^{\gamma_o} \backslash W/W^\gamma$. The next Lemma characterizes this subset in terms of the 
intersections $(\B[w])^\sigma=\B[w]\cap\Fcal^\R_\gamma$. 

\begin{lemma}\label{lem:schubert-reel} 
Let $w\in W$. The following statements are equivalent
\begin{enumerate}
\item $(\B[w])^{\sigma}\neq\emptyset$.
\item $(\Pbb[w])^{\sigma}\neq\emptyset$.
\item the class of $w$ in $W^{\gamma_o}\setminus W/W^\gamma$ is contained in $W^\sigma/W^\sigma\cap W^{\gamma}\simeq W_{\agot}/W_{\agot}^\gamma$.
\end{enumerate}
Moreover, if one of the statement is true, then $\B[w]=\Pbb[w]$ and $(\B[w])^\sigma= P[w]$.
\end{lemma}

{\em Proof :} The implication {\em 1.} $\Longrightarrow$ {\em 2.} is obvious since $\B[w]\subset \Pbb[w]$. For the implication \break 
{\em 2.} $\Longrightarrow$ {\em 3.} we consider 
the Bruhat decomposition $\Fcal_\gamma^\R=\bigcup_{v\in W_{\agot}/W_{\agot}^\gamma}P[v]$. As we have already explained, each class in 
$W_{\agot}/W_{\agot}^\gamma$ can be represented by an element of $W^\sigma$. If $(\Pbb[w])^\sigma\neq \emptyset$, there exists $v\in W^\sigma$ such that $P[v]\subset (\Pbb[w])^\sigma$. Hence $\Pbb[v]=\Pbb[w]$, and then $w=v$ in $W^{\gamma_o} \backslash W/W^{\gamma}$. The implication {\em 2.} $\Longrightarrow$ {\em 3.} is settled and the last one {\em 3.} $\Longrightarrow$ {\em 1.} is immediate.

Consider now $u\in W^\sigma$. Since $\B\subset\Pbb$, the Bruhat decomposition 1. tells us that $\Pbb[u]=\cup_{w}\B[w]$ where the union runs over the 
$w\in W/W^\gamma$ such that $w=u$ in $W^{\gamma_o}\setminus W/W^\gamma$. As $\sigma(u)=u$, the last relation implies that 
$w=u$ in $W/W^\gamma$. Thus we have proven that $\Pbb[u]=\B[u]$. Now the Bruhat decomposition 3. shows that 
$(\B[u])^\sigma=\Pbb[u]\cap \Fcal^\R_\gamma=\cup_{v} P[v]$ where the union runs over the $v\in W_{\agot}/W_{\agot}^\gamma\simeq W^\sigma/W^\sigma\cap W^{\gamma}$ such that  $P[v]\subset \Pbb[u]$. The last inclusion implies that $v=u$ in $W^{\gamma_o} \backslash W/W^{\gamma}$, and then 
$v=u$ in $W^\sigma/W^\sigma\cap W^{\gamma}$. We have proven that $(\B[u])^\sigma= P[u]$.
$\Box$

\subsection{Comparison with a result of Kapovich-Leeb-Millson}\label{sec:KLM}

Suppose that $\tG=G^s$ for $s\geq 2$ : Assumption \ref{hypothese} holds automatically. Here $\tpgot=\pgot^s$, $\pi(X_1,\cdots,X_s)=\sum_{j=1}^s X_j$, 
and the set $\Sigma(\tggot/\ggot)$ corresponds to the set $\Sigma$  of (non-zero) weights relative to the $\agot$-action on $\ggot$. The restricted Weyl group 
$W_{\tagot}$ is equal to $W_{\agot}^s$. 

Recall that we can associate to any $w\in W_{\agot}/W_\agot^\gamma$, the Schubert cell $\Xgot^o_{w,\gamma}:= \B[w]\subset \Fcal_\gamma$ and the Schubert variety $\Xgot_{w,\gamma}:=\overline{\Xgot^o_{w,\gamma}}$. 

In this setting, Theorem \ref{theo:main} becomes 

\begin{theorem}\label{theo:delta-exemple-KLM}  Let $(\xi_0,\xi_1,\cdots,\xi_s)\in(\agot_+)^{s+1}$. 
We have  $K\xi_0\subset\sum_{j=1}^s K\xi_j$ if and only if 
$$
\sum_{j=1}^s ( \xi_j,w_j\gamma)\geq ( \xi_0,w_0w\gamma)
$$
for any rational antidominant element $\gamma\in\agot$ and any $(w,w_1,\cdots,w_s)\in  
(W_\agot/W_\agot^\gamma)^{s+1}$ satisfying the following properties: 
\begin{enumerate}
\item[a)] ${\rm Vect}\big(\Sigma\cap \gamma^\perp\big)={\rm Vect}\big(\Sigma\big)\cap \gamma^\perp$.
\item[b)] $[\Xgot_{w,\gamma}]\cdot [\Xgot_{w_1,\gamma}]\cdot\ldots\cdot[\Xgot_{w_s,\gamma}]= [pt]$ in $H^*(\Fcal_\gamma,\Z)$.
\item[c)] ${\rm Tr}_{w\gamma}(\ngot^{w\gamma>0})+\sum_{j=1}^s{\rm Tr}_{w_j\gamma}(\ngot^{w_j\gamma>0})=s\, {\rm Tr}_{\gamma}(\ggot^{\gamma>0})$.
\end{enumerate}
\end{theorem}

The submanifold of $\Fcal_\gamma$ fixed by the conjugate linear involution  is naturally diffeomorphic to the real flag $\Fcal_\gamma^\R:=G/P_\lambda$.
 Thanks to the Bruhat decomposition $\Fcal_\gamma^\R=\bigcup_{w\in W_{\agot}/W_{\agot}^\gamma} P[w]$, we know that the cycle classes $[\Xgot^\R_{w,\gamma}], w\in W_{\agot}/W_{\agot}^\gamma$ associated to the real Bruhat-Schubert varieties $\Xgot_{w,\gamma}^\R:=\overline{P[w]}$  defines a basis of the cohomology $H^*(\Fcal_\gamma^\R,\Z_2)$ with $\Z_2$-coefficients \cite{Takeuchi-65}. 
 
 We know from Lemma \ref{lem:schubert-reel} that the real Bruhat-Schubert varieties $\Xgot_{w,\gamma}^\R$ corresponds to the real part of the complex 
 Bruhat-Schubert varieties $\Xgot_{w,\gamma}$. In this context, a theorem of Borel and Haefliger \cite{Borel-Haefliger} tells us that 
 $$
 b)\qquad [\Xgot_{w,\gamma}]\cdot [\Xgot_{w_1,\gamma}]\cdot\ldots\cdot[\Xgot_{w_s,\gamma}]= [pt]\quad  in\quad 
 H^*(\Fcal_\gamma,\Z).
 $$ 
 implies the following
 $$
 b^{\R})\qquad [\Xgot^\R_{w,\gamma}]\cdot [\Xgot^\R_{w_1,\gamma}]\cdot\ldots\cdot[\Xgot^\R_{w_s,\gamma}]= [pt]\quad  in\quad 
 H^*(\Fcal_\gamma^\R,\Z_2).
 $$
 
The description of $\Pi_{\pgot}(K^s,K)$ obtained by Kapovich-Leeb-Millson \cite{KLM-memoir-08} was in terms of elements
$(\gamma,w,w_1,\cdots,w_s)\in\agot\times W_\agot^{s+1}$ verifying conditions $a)$ and $b^{\R})$. Our description is therefore more accurate, first by adding condition $c)$ and then by taking the refined condition $b)$. 

\subsection{Levi-movability}\label{sec:levi}

In this section, we recall Belkale-Kumar's notion of Levi-movability \cite{BK1}  and explain its connection with conditions $b)$ and $c)$ of Theorem \ref{theo:delta-exemple-KLM}.

Recall that we work with a reductive Lie group $U_\C$ equipped with an anti-linear involution $\sigma$. The Lie algebra $\ggot:=(\ugot_\C)^\sigma$ is a real form of the complex Lie algebra $\ugot_\C$. We have chosen a maximal torus $T\subset U$ with Lie algebra $\tgot$, that is invariant under $\sigma$, and such that $\tgot^{-\sigma}$ is of maximal dimension. Let $\agot=\frac{1}{i}\tgot^{-\sigma}\subset (\ugot_\C)^\sigma$.

We choose a Weyl chamber $\tgot_+\subset \tgot$ such that $\agot_+=\frac{1}{i}(\tgot^{-\sigma}\cap\tgot_+)$ is a Weyl chamber relatively to the restricted root system $\Sigma$. Let $\B\subset U_\C$ be the Borel subgroup with Lie algebra $\tgot_\C\oplus \sum_{\alpha \in \Rgot^+}(\ugot_\C)_\alpha$.
To any element $\gamma\in \agot$, we associate the parabolic subgroup $\Pbb_\gamma=\{g\in U_\C, \lim_{t\to\infty}e^{t\gamma}ge^{-t\gamma}\ {\rm exists}\}$ that is invariant under the involution $\sigma$. Notice that $\B\subset \Pbb_\gamma$ when $\gamma$ is anti-dominant, i.e. $\gamma\in -\agot_+$.

The antilinear involution $\sigma$ defines an anti-holomophic involution on the flag variety $\Fcal_\gamma= U_\C/\Pbb_\gamma$. We consider 
the cellular decomposition 
$\Fcal_\gamma=\amalg_{w\in W/W^\gamma} \B [w]$ parameterized by the Weyl group $W=N_U(\tgot)/T$ and the subgroup $W^\gamma$ fixing 
$\gamma\in \agot$. For any $w\in W/W^\gamma$, we denote by $\Xgot_{w,\gamma}^o=\B [w]\subset \Fcal_\gamma$ the Schubert cell.

The cycle class of the Schubert variety $\Xgot_{w,\gamma}=\overline{\B [w]}$ in $H^*(\Fcal_\gamma,\Z)$ is denoted by $[\Xgot_{w,\gamma}]$ and it is called a Schubert class. Let $\Xgot_{w_i,\gamma}$ be $d$ Schubert classes parametrized by $w_i \in W/W^\gamma$, $i\in [d]$.  
If there exists an integer $k$ such that $[\Xgot_{w_1,\gamma}]\cdot\ldots\cdot[\Xgot_{w_d,\gamma}]= k [pt]$, then we set $c_\gamma(w_1,\ldots,w_d) =k$; 
we set $c_\gamma(w_1,\ldots,w_d)=0$ otherwise.

Using the transversality theorem of Kleiman, Belkale and Kumar showed in \cite{BK1},  Proposition 2,  the following important lemma.
\begin{lemma}Suppose that $\gamma$ is anti-dominant.
The coefficient $c_\gamma(w_1,\ldots,w_d)$ is nonzero if and only if for generic $(p_1,\ldots,p_d)\in (\Pbb_\gamma)^d$, the intersection
$p_1w_1^{-1}\Xgot_{w_1,\gamma}^o \cap \cdots \cap p_dw_d^{-1}\Xgot_{w_d,\gamma}^o$
is transverse at $e$. 
\end{lemma}

Then Belkale-Kumar defined Levi-movability.

\begin{definition}Suppose that $\gamma$ is anti-dominant. A d-uplet $(w_1,\ldots,w_d) \in W^d$ is Levi-movable if for generic $(l_1,\ldots,l_d)\in (U^\gamma_\C)^d$, the intersection
$l_1w_1^{-1}\Xgot_{w_1,\gamma}^o \cap \cdots \cap l_dw_d^{-1}\Xgot_{w_d,\gamma}^o$
is transverse at $e$. 
\end{definition}

Belkale-Kumar have defined in \cite{BK1}, a new product $\odot_0$ on $H^*(\Fcal_\gamma,\Z)$  that is closely related to the notion of Levi-movability. By definition of this product, the identity
$$
[\Xgot_{w_1,\gamma}]\odot_0\ldots\odot_0[\Xgot_{w_d,\gamma}]= \ell [pt] 
$$ 
holds for some $\ell\geq 1$ if and only if 
$[\Xgot_{w_1,\gamma}]\cdot\ldots\cdot[\Xgot_{w_d,\gamma}]= \ell [pt]$ and the $d$-uplet $(w_1,\ldots,w_d)$ is Levi-movable.

In the next proposition we concentrate ourself to $d$-uplets $(w_1,\ldots,w_d)$ of $\sigma$-invariant elements. Recall that $\ngot\subset \ggot$ is the nilpotent Lie subalgebra defined as the subspace $\ggot^{\gamma_o>0}$ for any regular dominant element $\gamma_o\in \agot_+$.

\begin{proposition}\label{prop:condition-b-c}Let $\gamma\in\agot$ be a rational anti-dominant element.
A d-uplet $(w_1,\ldots,w_d)\in (W^\sigma)^d$ satisfies the following conditions 
\begin{enumerate}
\item[b)] $[\Xgot_{w_1,\gamma}]\cdot\ldots\cdot[\Xgot_{w_d,\gamma}]= [pt]$ in $H^*(\Fcal_\gamma,\Z)$.
\item[c)] $\sum_{j=1}^d{\rm Tr}_{w_j\gamma}(\ngot^{w_j\gamma>0})=(d-1) {\rm Tr}_{\gamma}(\ggot^{\gamma>0})$.
\end{enumerate}
if and only if $[\Xgot_{w_1,\gamma}]\odot_0\ldots\odot_0[\Xgot_{w_d,\gamma}]= [pt]$.
\end{proposition}

{\em Proof:} Since $\ggot$ is the real part of the  complex Lie algebra $\ugot_\C$, we have the following identity of $\agot$-modules : 
$\ugot_\C\simeq \ggot^2$. Thus ${\rm Tr}_{\gamma}(\ggot^{\gamma>0})=\tfrac{1}{2}{\rm Tr}_{\gamma}(\ugot_\C^{\gamma>0})$. 
The nilpotent radical of the Lie algebra of $\B$ is $\ngot_{\ugot_\C}= E\oplus \ngot_\C$, where $E=\ngot_{\ugot_\C}\cap Z_{\ugot_\C}(\agot)$. Thus for any $w\in W_\agot$, the vector space $\ngot_{\ugot_\C}^{w\gamma>0}$ is 
the complexification of $\ngot^{w\gamma>0}$. It follows that ${\rm Tr}_{w\gamma}(\ngot^{w\gamma>0})=
\tfrac{1}{2}{\rm Tr}_{w\gamma}(\ngot_{\ugot_\C}^{w\gamma>0})$. We have proven that $c)$ is equivalent to 
$\sum_{j=1}^d{\rm Tr}_{w_j\gamma}(\ngot_{\ugot_\C}^{w_j\gamma>0})=(d-1){\rm Tr}_{\gamma}(\ugot_\C^{\gamma>0})$. 
Let $\rho_\ugot$ be half the sum of the positive roots and let $\gamma'=i\gamma\in\tgot^{-\sigma}\cap -\tgot_+$. Standard computations give that 
${\rm Tr}_{w\gamma}(\ngot_{\ugot_\C}^{w\gamma>0})=\langle\rho_\ugot,w\gamma'\rangle -\langle\rho_\ugot,\gamma'\rangle$ and 
${\rm Tr}_{\gamma}(\ugot_\C^{\gamma>0})=-2\langle\rho_\ugot,\gamma'\rangle$. Finally we see that $c)$ is equivalent to 
$$
c')\qquad\sum_{j=1}^d\langle\rho_\ugot,w_j\gamma'\rangle=(2-d)\langle\rho_\ugot,\gamma'\rangle.
$$

Consider the complex $U_\C$-manifold $M=(U_\C/\B)^d$. The set 
$C=\Pi_{k=1}^d U^{\gamma'}_\C[w_k^{-1}]$ is a connected component of the submanifold $M^{\gamma'}$. The Bialynicki-Birula's complex submanifold 
$C^-:=\{m\in M,  \lim_{t\to\infty} \exp(-it\gamma') m\ \in C\}$ is equal to $C=\Pi_{k=1}^d \Pbb'[w_k^{-1}]$ where $\Pbb'$ is the parabolic subgroup with Lie algebra $\tgot_\C\oplus \sum_{\langle\alpha,\gamma'\rangle \leq 0}(\ugot_\C)_{\alpha}$. Notice that the Borel subgroup $\B$ is contained in $\Pbb'$ as 
$\gamma'\in -\tgot_+$.

Now we consider the map $\pi:U_\C\times_{\Pbb'} C^-\to M$ that sends $[g,x]$ to $gx$. Standard computations shows the following statement are equivalent 
\begin{enumerate}
\item Conditions $b)$ and $c')$ holds.
\item $\pi$ is a birationnal map with an exceptional set that does not contains $C$.
\item $[\Xgot_{w_1,\gamma}]\odot_0\ldots\odot_0[\Xgot_{w_d,\gamma}]= [pt]$.
\end{enumerate}
See \S 5 in \cite{Brion-Bourbaki}. $\Box$

\medskip

Proposition \ref{prop:condition-b-c} permits us to describe the convex cone $\Pi_{\pgot}(K^s,K)$ in terms of the Belkale-Kumar's product $\odot_0$.

\begin{theorem}\label{theo:version-levi-movable}   
Let $(\xi_0,\xi_1,\cdots,\xi_s)\in(\agot_+)^{s+1}$. 
We have  $K\xi_0\subset\sum_{j=1}^s K\xi_j$ if and only if  $\sum_{j=1}^s ( \xi_j,w_j\gamma)\geq ( \xi_0,w_0w\gamma)$
for any rational antidominant element $\gamma\in\agot$ and any $(w,w_1,\cdots,w_s)\in  
(W_\agot/W_\agot^\gamma)^{s+1}$ satisfying the following properties: 
\begin{enumerate}
\item  ${\rm Vect}\big(\Sigma\cap \gamma^\perp\big)={\rm Vect}\big(\Sigma\big)\cap \gamma^\perp$.
\item $[\Xgot_{w,\gamma}]\odot_0[\Xgot_{w_1,\gamma}]\odot_0\ldots\odot_0[\Xgot_{w_s,\gamma}]= [pt]$ in $H^*(\Fcal_\gamma,\Z)$.
\end{enumerate}
\end{theorem}

\subsection{Singular Horn cone}\label{sec:singular-horn}

Let $n\geq 1$. If $A$ is an Hermitian $n\times n$ matrix, we denote by $\mathrm{s}(A)=(\mathrm{s}_1(A)\geq \cdots \geq\mathrm{s}_n(A))$ its spectrum. The Horn cone 
$\horn(n)$ is defined as the set of triplets $(\mathrm{s}(A),\mathrm{s}(B),\mathrm{s}(C))$ where $A,B,C$ are Hermitian $n\times n$ matrices satisfying 
$A+B+C=0$.

Denote the set of cardinality $r$ subsets $I=\{i_1<i_2<\cdots<i_r\}$ of $[n]=\{1,\ldots,n\}$ by $\Pcal^n_r$. To each $I\in \Pcal^n_r$
we associate 
\begin{itemize}
\item a weakly decreasing sequence of non-negative integers 
$\lambda(I)=(\lambda_1\geq\cdots\geq \lambda_r)$ where $\lambda_a= n-r+a-i_a$ for $a\in [r]$.
\item a Schubert class $[\Xgot_I]\in H^*(\G(r,n),\Z)$ where $\G(r,n)$ denotes the Grassmann variety of $r$-dimensional subspaces of $\C^n$.
\end{itemize}

\medskip

The following Horn's conjecture \cite{Horn} was settled in the affirmative by combining the work of A. Klyachko \cite{Klyachko} with the work of A. Knutson and T. Tao 
\cite{Knutson-Tao-99} on the ``saturation'' problem. We refer the reader to survey articles \cite{Fulton00,Brion-Bourbaki} for details.

If $x=(x_1,\ldots,x_n)\in\R^n$ and $I\subset [n]$, we define $|\,x\,|_I=\sum_{i\in I}x_i$ and $|\,x\,|=\sum_{i=1}^n x_i$. Let $1^r=(1,\ldots,1)\in \R^r$. 

\begin{theorem}[Horn's conjecture]\label{theo:horn-conjecture}
An element $(x,y,z)\in(\R^n_+)^3$ belongs to $\horn(n)$ if and only if the following conditions holds : $|\,x\,|+|\, y\,|+|\,z\,|=0$ 	and 
$$
|\,x\,|_{I}+|\,y\,|_{J} +|\, z\,|_{K}\leq 0
$$
for any $r\in [n-1]$, for any $I,J,K\in\Pcal^n_r$ such that 
\begin{equation}\label{eq:condition-horn-r}
(\lambda(I),\lambda(J),\lambda(K)-(n-r) 1^r)\in\horn(r).
\end{equation}
\end{theorem}

Thanks to the saturation Theorem of Knutson-Tao, we know that (\ref{eq:condition-horn-r}) is equivalent to the following relation in $H^*(\G(r,n),\Z)$ 
\begin{equation}\label{eq:condition-produit}
[\Xgot_I]\cdot[\Xgot_J]\cdot[\Xgot_K]=\ell [pt]\quad \ell\geq 1,
\end{equation}
where $[pt] \in H^{{\small\rm max}}(\G(r,n),\Z)$ denote the Poincar\'e dual class of the point. Belkale proved in \cite{Belkale06} that Theorem 
\ref{theo:horn-conjecture} holds if we replace condition (\ref{eq:condition-produit}) by
\begin{equation}\label{eq:condition-produit=1}
[\Xgot_I]\cdot[\Xgot_J]\cdot[\Xgot_K]= [pt].
\end{equation}

Finally, Knutson-Tao-Woodward proved in \cite{KTW04} that the list of inequalities parameterized by the $(r,I,J,K)$ satisfying (\ref{eq:condition-produit=1}) is optimal.

\medskip

We turn now our attention to the singular Horn cone. Suppose that $n=p+q$ with $p\geq q\geq 1$.  If $A$ is a complex $p\times q$ matrix, we denote by 
$\tau(A)=(\tau_1(A)\geq \cdots \geq\tau_q(A)\geq 0)$ its singular spectrum : $\tau_k(A)=\sqrt{\mathrm{s}_k(A^*A)}$ for any $k\in [q]$.

\begin{definition}
The singular Horn cone, denoted  $\sing(p,q)$, is defined as the set of triplets $(\tau(A),\tau(B),\tau(C))$ where $A,B,C$ are complex 
$p\times q$ matrices satisfying $A+B+C=0$.
\end{definition}

We have a natural connection between $\sing(p,q)$ and $\horn(n)$ because, for any complex $p\times q$ matrix A, the spectrum of the $n\times n$ matrix 
$$
\begin{pmatrix} 
0 & A\\
A^* &0
\end{pmatrix}
$$
is $(\tau_1(A)\geq \cdots \geq\tau_q(A)\geq 0\geq\cdots\geq 0\geq -\tau_q(A)\geq \cdots \geq -\tau_1(A))$. To any $x=(x_1,\ldots, x_q)\in \R^q$ we associate 
the vector $\tilde{x}=(x_1,\ldots, x_q,0,\ldots,0,- x_q,\ldots,- x_1)\in\R^n$. From the previous remark, we see that 
$$
(x,y,z)\in \sing(p,q)\Longrightarrow (\tilde{x},\tilde{y},\tilde{z})\in \horn(n).
$$
The O'Shea-Sjamaar theorem tells us that the implication above is in fact an equivalence. Hence we obtain a first description of $\sing(p,q)$.

\begin{proposition} 
Let $\R^q_{++}=\{x=(x_1\geq \cdots\geq x_q\geq 0)\}$. An element $(x,y,z)\in(\R^q_{++})^3$ belongs to $\sing(p,q)$ if and only if 
\begin{equation}\label{eq:condition-singular}
|\,\tilde{x}\,|_{I}+|\,\tilde{y}\,|_{J} +|\, \tilde{z}\,|_{K}\leq 0
\end{equation}
for any $r\in [n-1]$, for any $I,J,K\in\Pcal^n_r$ such that $[\Xgot_I]\cdot[\Xgot_J]\cdot[\Xgot_K]= [pt]$ in $H^*(\G(r,n),\Z)$.
\end{proposition}

The main issue with this first description is that most of the inequalities (\ref{eq:condition-singular}) are redundant. 

Now we explain the more precise description of $\sing(p,q)$ that we obtain by applying Theorem \ref{theo:version-levi-movable}. For any $r\in [q]$, let 
$\Pcal^{p,q}_r\subset \Pcal^n_r$ be the collection of subsets $I$ satisfying\footnote{$I^o=\{n+1-k, k\in I\}$} $I\cap I^o=\emptyset$ and $I\cap \{q+1,\ldots, p\}=\emptyset$. 
We see that any $I\in \Pcal^{p,q}_r$ is equal to the union $I_+\coprod I_{-}^o$ where $I_{+},I_{-}$ are disjoint subsets of $[q]$.

Notice that, when $I,J,K\in\Pcal^{p,q}_r$, the inequality (\ref{eq:condition-singular}) becomes
$$
(\star)_{I_\pm,J_\pm,K_\pm}\hspace{2cm}|\,x\,|_{I_+}+|\,y\,|_{J_+} +|\, z\,|_{K_+}\leq |\,x\,|_{I_{-}}+|\,y\,|_{J_{-}} +|\, z\,|_{K_{-}}.
$$

The first version of our result is the following theorem that answers (partially) a conjecture of A. S. Buch (see \cite{Fulton00}, \S 5).
\begin{theorem}
An element $(x,y,z)\in(\R^q_{++})^3$ belongs to $\sing(p,q)$ if and only if $(\star)_{I_\pm,J_\pm,K_\pm}$ holds 
for any $r\in [q]$, for any $I,J,K\in\Pcal^{p,q}_r$ such that $[\Xgot_I]\cdot[\Xgot_J]\cdot[\Xgot_K]= [pt]$ in $H^*(\G(r,n),\Z)$.
\end{theorem}

The above description can be refined as follows. For any $r\in [q]$ such that $r<\frac{n}{2}$, we denote by $\F(r,n-r;n)$ the two-steps flag variety  parameterizing 
nested sequences of  linear subspaces $E\subset F^\perp\subset \C^n$ where $\dim E=\dim F= r$. Here the orthogonal is taken relatively to the bilinear product 
$(x,y)=\sum_{i=1}^n x_iy_i$. In the case where $p=q$ and $r=q=\frac{n}{2}$, we denote by $\F(q,q;n)$  (or by $\G(q,n)$) the Grassmannian  parameterizing linear subspaces $E\subset \C^n$ of dimension $q$.

\begin{definition}
Let $A\subset [n]$ be a subset of cardinal $2r$. We say that $A$ is polarized if it admits a decomposition 
$A=A_+ \coprod A_-$ into disjoints subsets. If  $A_+$ and $A_-$ are both of cardinal $r$, we say that $A$ is balanced. We denote by $\Pol^n_{2r}$ the collection of balanced polarized subsets of $[n]$ of cardinal $2r$.
\end{definition}

The orbits of the Borel subgroup $B_n\subset GL_n(\C)$ of upper-triangular matrices on $\F(r,n-r;n)$ are parametrized by the balanced polarized subsets of $[n]$ of cardinal $2r$. To any balanced polarized subset $A=A_+ \coprod A_-\in \Pol^n_{2r}$, we associate  
\begin{itemize}
\item the flag $\C^{A_+}\subset (\C^{A_-})^\perp$, where $\C^{A_\pm}={\rm Vect}(e_i, i\in A_\pm)\subset \C^n$. Notice that in the case where $p=q$ and $r=q=\frac{n}{2}$, we have $\C^{A_+}=(\C^{A_-})^\perp$.
\item the Schubert variety $\Xgot_A=\overline{B_n\cdot(\C^{A_+}\subset (\C^{A_-})^\perp)}\subset \F(r,n-r;n)$ , and 
its cycle class $[\Xgot_A]\in H^*(\F(r,n-r;n),\Z)$ in cohomology. 
\end{itemize}

To any subset $I\in\Pcal^{p,q}_r$, we associate the balanced polarized subset $\widehat{I}=I\coprod I^o\in  \Pol^n_{2r}$ and the cohomology class 
$[\Xgot_{\widehat{I}}]\in H^*(\F(r,n-r;n),\Z)$. The final description of  $\sing(p,q)$ is given in the following Theorem (whose proof is postponed at  
\S \ref{sec:preuve-singular-horn}).

\begin{theorem}\label{theo:singular-p-q}
An element $(x,y,z)\in(\R^q_{++})^3$ belongs to $\sing(p,q)$ if and only if 
$(\star)_{I_\pm,J_\pm,K_\pm}$ holds for any $r\in [q]$, for any $I,J,K\in\Pcal^{p,q}_r$ such that 
\begin{equation}\label{eq:condition-cohomologique-2-flag-section}
[\Xgot_{\widehat{I}}]\odot_0[\Xgot_{\widehat{J}}]\odot_0[\Xgot_{\widehat{K}}]= [pt]\quad {\rm in}\quad  H^*(\F(r,n-r;n),\Z).
\end{equation}
\end{theorem}

\medskip

In Section \S \ref{sec:multiplicative}, we will explain why condition (\ref{eq:condition-cohomologique-2-flag-section}) implies condition (\ref{eq:condition-produit=1}).

\medskip

\medskip 

Let us end this section with the example of $\sing(3,3)$. The first description obtained from the O'Shea-Sjamaar theorem says that 
$\sing(3,3)$ admits a natural identification with the intersection $\horn(6)\cap V^3$ where $V=\{(a,b,c,-c,-b,-a)\}\subset \R^6$. Thanks to computer calculation we know that $\horn(6)\subset (\R^6)^3$ is described with a minimal list of 536 inequalities. Theorem \ref{theo:singular-p-q} allows to describe $\sing(3,3)$ in terms of 87 $+$ 9 inequalities (see \S \ref{sec:example-horn-1-2-3}).

\section{Parametrization of facets of real Kirwan polytopes}

Let $L_\C$ be a connected complex reductive group with maximal compact subgroup $L$.  
Let $\sigma$ be a conjugate linear involution on $L_\C$ commuting with the Cartan involution.
We work with a $\sigma$-invariant scalar product $(-,-)$ on $\lgot$, allowing an identification $\lgot\simeq \lgot^*$.

Let us consider a proper K\"{a}hler-Hamiltonian $L$-manifold with involution $(M,\Omega,\J,\Phi_\lgot,\tau,\sigma)$ (see \S \ref{sec:kahler-hamiltonian}).  By that we mean~:
\begin{itemize}
\item $L_\C$ acts holomorphically on the complex manifold $(M,\J)$ and the K\"{a}hler form $\Omega$ is $L$-invariant, 
\item There is a $L$-equivariant {\em proper} moment mapping $\Phi_\lgot : M\to \lgot$ satisfying the relation\footnote{Here $X_M(m)=-X\cdot m =\frac{d}{dt}\vert_{t=0}e^{-tX}m$ is the vector field generated by $X\in \lgot$.} $d(\Phi_\lgot,X)=-\Omega(X_M, -)$,  $\forall X\in\lgot$.
\item The anti-holomorphic involution $\tau :M\to M$ satisfies the following relations : $\tau^*(\Omega)=-\Omega$, and $\tau(g\cdot m)=\sigma(g)\cdot\tau(m)$, 
$\Phi_\lgot(\tau(m))=-\sigma(\Phi_\lgot(m))$, for all $(g,m)\in L_\C\times M$. 
\end{itemize}

Let us denote by $G_L$ the connected component of the subgroup $(L_\C)^\sigma$. The Lie algebra of $G_L$ admits the Cartan decomposition $\ggot_L=\kgot_L\oplus \pgot_L$ where $\pgot_L=i\lgot^{-\sigma}$ and $\kgot_L$ is the Lie algebra of $K_L:=(L^\sigma)_o$. 

The submanifold $Z=\{m\in M, \tau(m)=m\}$, that we suppose {\em non-empty}, is a Lagrangian submanifold of $(M,\Omega)$, not necessarily connected, and stable under the action of the real reductive group $G_L$. Let $\Zcal$ be a connected component of $Z$. We define the $K_L$-equivariant map 
$$
\Psi : \Zcal\longrightarrow \pgot_L
$$
by the relation $\Psi(z)=\frac{1}{i}\Phi_\lgot(z)$. One sees that, for any $\beta\in\pgot_L$,  $\beta_\Zcal$ is the gradient vector field of the function $(\Psi,\beta)$ .

 The O'Shea-Sjamaar Theorem \cite{OSS99} tells us that the image $\Psi(\Zcal)$ is equal to $\Phi_\lgot(M)\bigcap \lgot^{-\sigma}$ through the identification 
 $\pgot_L\simeq \lgot^{-\sigma}$ (for another proof see \cite{pep-real-RP1}). Let $\agot_L\subset \pgot_L$ be a maximal abelian subspace. 
 The $K_L$-orbits of $\pgot_L$ are parametrized by a Weyl chamber $(\agot_L)_+\subset\agot_L$. Hence 
the $K_L$-orbits of $\Psi(\Zcal)$ are parametrized by the real Kirwan polytope
$$
\Delta(\Zcal):=\Psi(\Zcal)\bigcap \, (\agot_L)_+.
$$
Thanks to the O'Shea-Sjamaar Theorem, we know that $\Delta(\Zcal)$ is a convex polytope. In the next sections we explain the parametrization of the facets 
of $\Delta(\Zcal)$ in terms of balanced Ressayre's pairs that we obtained in \cite{pep-real-RP1}.

\subsection{Admissible elements} 

The stabilizer subgroups of $m\in M$ relatively to the $L$ and $G_L$ actions are denoted respectively by $L^m$ and $(G_L)^m$ : their Lie algebras are denote by 
$\lgot^m$ and $(\ggot_L)^m$. We will take a particular attention to the subspace $(\pgot_L)^m=\{X\in\pgot_L, X\cdot m=0\}\subset\ggot_m$. 
 If $\gamma\in\agot_L$, we denote by $\Zcal^\gamma$ the submanifold where the vector field $z\mapsto \gamma\cdot z$ vanishes.

\begin{definition} Let us define 
$$
\dim_{\pgot_L}(\Xcal):=\min_{z\in\Xgot}\,\dim((\pgot_L)^z)
$$
for any subset $\Xcal\subset \Zcal$. A non-zero element $\gamma\in\agot$ is called {\em admissible} relatively to the $G$-manifold $\Zcal$  if $\gamma$ is rational, and if 
$\dim_{\pgot_L}(\Zcal^\gamma)-\dim_{\pgot_L}(\Zcal)\in\{0,1\}$.
\end{definition}

The next results, proven in \cite{pep-real-RP1}, are useful to determine  $\dim_{\pgot_L}(\Zcal)$.

\begin{lemma}\label{lem:stabilisateur-min}
There exists a subspace $\hgot\subset \pgot_L$ such that 
\begin{enumerate}
\item $\forall x\in \Zcal, \exists k\in K_L$ such that ${\rm Ad}(k)(\hgot)\subset (\pgot_L)^x$,
\item $\dim(\hgot)=\dim((\pgot_L)^x)$ on a dense open subset of $\Zcal$.
\end{enumerate}
Therefore $\dim_{\pgot_L}(\Zcal)=\dim(\hgot)$.
\end{lemma}

The subspace $\hgot$ is called {\em the generic $\pgot_L$-stabilizer} of $\Zcal$.

\begin{corollary}\label{coro:stabilisateur-min}
Let $\hgot_o\subset \pgot_L$ be a subspace such that the set $K \Zcal^{\hgot_o}$ has a non-empty interior in $\Zcal$. Then 
$\dim_{\pgot_L}(\Zcal)\geq \dim(\hgot_o)$.
\end{corollary}

\subsection{Balanced Ressayre's pair}\label{sec:balanced-RP}

Let $\gamma\in\agot_L$ be a non-zero element, and let $C\subset M^\gamma$ be a connected component  intersecting $\Zcal$. The K\"{a}hler submanifold $C$ is stable under the involution $\tau$ and under the action of the stabilizer sugbroup $L^\gamma_\C=\{g\in L_\C, g\gamma=\gamma\}$. 
The real submanifold $\Ccal= C\cap \Zcal\subset \Zcal^\gamma$ 
is non-empty, not necessarily connected, and stable under the action of the group $G_L^\gamma=L^\gamma_\C\cap G_L$. The fonction 
$(\Psi,\gamma)$ is locally constant on $M^\gamma$  and we denote by $(\Psi(\Ccal),\gamma)$ its value on $\Ccal$.

The Bialynicki-Birula's complex submanifold 
\begin{equation}\label{eq:BB}
C^-:=\{m\in M,  \lim_{t\to\infty} \exp(t\gamma) m\ \in C\}.
\end{equation}
is invariant under the involution of $M$. The infinitesimal action of $\gamma$ defines a symmetric endomorphism of the tangent bundle $\T M\vert_C$, and one sees that 
$\T C^-\vert_C= (\T M\vert_C)^{\gamma\leq 0}$.

Consider now the parabolic subgroups defined by
\begin{equation}\label{eq:P-gamma-reel}
\Pbb_\gamma=\{g\in L_\C, \lim_{t\to\infty}\exp(t\gamma)g\exp(-t\gamma)\ {\rm exists}\}\quad {\rm and}\quad P_\gamma=G_L\cap \Pbb_\gamma.
\end{equation}
Here the Lie algebra of $\Pbb_\gamma$ and $P_\gamma$ are respectively $(\lgot_\C)^{\gamma\leq 0}$ and $(\ggot_L)^{\gamma\leq 0}$. We denote by $\Pbb\subset  L_\C$ and $P\subset G_L$ the parabolic subgroups $\Pbb_{-\gamma_o}$ and $P_{-\gamma_o}$, where $\gamma_o$ is an element of the 
{\em interior} of the Weyl chamber $(\agot_L)_+$. The nilradical of the Lie algebra of $P$ is $\ngot_L:=\ggot^{\gamma_o>0}$.

The submanifold $\Ccal^-$ is invariant under the action of $\Pbb_\gamma$, hence we can consider 
the complex manifold $\Pbb\times_{\Pbb\cap \Pbb_\gamma} C^-$ and the holomorphic map 
$$
{\rm q}^\C_\gamma: \Pbb\times_{\Pbb\cap \Pbb_\gamma} C^- \to M
$$
that sends $[p,m]$ to $pm$. The complex manifold $\Pbb\times_{\Pbb\cap \Pbb_\gamma} C^-$ is equipped with a natural anti-holomorphic involution 
$\tau_\gamma:[p,m]\mapsto [\sigma(p),\tau(m)]$ such that $\tau\circ {\rm q}^\C_\gamma ={\rm q}_\gamma^\C\circ \tau_\gamma$. The submanifold fixed by the involution is  $P\times_{P\cap P_\gamma} \Ccal^-$ where $\Ccal^-=C^-\cap \Zcal$ is a submanifold of $\Zcal$. Notice that for any $z\in \Ccal$, we have 
$\T_z \Ccal^-=(\T_z \Zcal)^{\gamma\leq 0}$.

We can now define the {\em real} counterpart of the map ${\rm q}^\C_\gamma$ : the map
$$
{\rm q}_\gamma: P\times_{P\cap P_\gamma} \Ccal^- \to \Zcal.
$$
The tangent bundle at the identity of $P/P\cap P_\gamma$ admits a canonical identification with $(\ngot_L)^{\gamma>0}$, so that for any $z\in \Ccal$, 
the tangent map $\T{\rm q}_\gamma\vert_{[e,z]}:(\ngot_L)^{\gamma>0}\times (\T_z \Zcal)^{\gamma\leq 0}\longmapsto \T_z \Zcal$ 
 is defined by the relation $\T{\rm q}_\gamma\vert_{[e,z]}(X,v)=X\cdot z\oplus v$. Hence for any $z\in\Ccal$ the following statements are equivalent :
 \begin{itemize}
  \item $\T{\rm q}^\C_\gamma\vert_{[e,z]}$ is an isomorphism
 \item $\T{\rm q}_\gamma\vert_{[e,z]}$ is an isomorphism
 \item $\rho_z^\gamma :(\ngot_L)^{\gamma>0}\longmapsto (\T_z \Zcal)^{\gamma> 0}$, $X\mapsto X\cdot z$ is an isomorphism.
 \end{itemize}

\begin{definition}\label{def:ressayre-pair}
The data $(\gamma,\Ccal)$ is called a {\bf balanced Ressayre's pair} of $\Zcal$ if there exists $\tau$-invariant  Zariski open subsets 
$\Vcal^\C\subset M$ and $\Ucal^\C\subset C^-$ such that : $\Vcal^\C$ is  $\Pbb$-invariant, \break $\Ucal^\C$ is  $\Pbb\cap \Pbb_{\gamma}$-invariant and intersects $C$, 
and the map ${\rm q}^\C_{\gamma}$ defines a diffeomorphism \break $\Pbb\times_{\Pbb\cap \Pbb_{\gamma}}\Ucal^\C\simeq \Vcal^\C$.

If furthermore we have  $\dim_\pgot(\Ccal)-\dim_\pgot(\Zcal)\in\{0,1\}$, and $\gamma$ is rational, we call $(\gamma,\Ccal)$ a regular balanced Ressayre's pair.
\end{definition}

If $\Ucal^\C$ is a $\tau$-invariant Zariski open subset  of $C^-$ that intersects the complex submanifold $C$, then 
$\Ucal^\C\cap C$ is a $\tau$-invariant Zariski open subset of $C$. It implies that $\Ucal^\C\cap C^\tau$ is dense in $C^\tau$, and in particular that 
$\Ucal^\C\cap \Ccal$ is non-empty (see Appendix B in \cite{pep-real-RP1}).  When the map ${\rm q}^\C_{\gamma}$ defines a diffeomorphism 
$\Pbb\times_{\Pbb\cap \Pbb_{\gamma}}\Ucal^\C\simeq \Vcal^\C$, we see then that the infinitesimal map  $ \rho_z^\gamma :(\ngot_L)^{\gamma>0}\longmapsto (\T_z \Zcal)^{\gamma> 0}$ is an isomorphism for any $z\in \Ucal^\C\cap \Ccal$.

We can weaken the condition on the balanced Ressayre's pair.

\begin{definition}\label{def:ressayre-pair-infinitesimal}
The data $(\gamma,\Ccal)$ is called a {\bf infinitesimal Ressayre's pair} of $\Zcal$ if there exists $z\in\Ccal$ such that $\rho_z^\gamma :(\ngot_L)^{\gamma>0}\longmapsto (\T_z \Zcal)^{\gamma> 0}$ is an isomorphism. If furthermore we have  $\dim_\pgot(\Ccal)-\dim_\pgot(\Zcal)\in\{0,1\}$, and 
$\gamma$ is rational, we call $(\gamma,\Ccal)$ a regular infinitesimal Ressayre's pair of $\Zcal$.
\end{definition}

When $(\gamma, \Ccal)$ is an infinitesimal real Ressayre's pair of $\Zcal$, there exists $z\in\Ccal$ so that $\rho_z^\gamma:(\ngot_L)^{\gamma>0}\longrightarrow  (\T_x \Zcal)^{\gamma>0}$ is an isomorphism commuting with the infinitesimal action of $\gamma$. It gives us two identities 
\begin{eqnarray*}
({\rm C}_1)&\qquad& \qquad\dim (\ngot_L)^{\gamma>0}={\rm rank}(\T\Zcal\vert_\Ccal)^{\gamma>0},\\
({\rm C}_2)&\qquad&\qquad{\rm Tr}_\gamma ((\ngot_L)^{\gamma>0})={\rm Tr}_\gamma((\T\Zcal\vert_\Ccal)^{\gamma>0}).
\end{eqnarray*}

If $(\gamma,\Ccal)$ is a balanced Ressayre's pair then 
$$
({\rm C}_3) \qquad \left\{x\in M, \sharp \{({\rm q}_\gamma^\C)^{-1}(x)\}=1\right\} \ {\rm contains\ a\ Zariski\ open\ subset\ of}\ M.
$$

Condition $({\rm C}_3)$ can be made more flexible by considering
$$
({\rm C}'_3) \qquad \left\{x\in M, ({\rm q}_\gamma^\C)^{-1}(x)\neq \emptyset \right\} \ {\rm contains\ a\ Zariski\ open\ subset\ of}\ M.
$$

\subsection{Parametrization of the facets}

The following result will be the main tool used in the next sections.

\begin{theorem}[\cite{pep-real-RP1}]\label{theo-facet-kirwan-reel}
For $\xi\in(\agot_L)_{+}$, the following statements are equivalent:
\begin{itemize}
\item $\xi\in\Delta_\pgot(\Zcal)$.
\item $(\xi,\gamma)\geq (\Psi(\Ccal),\gamma)$ holds for any regular infinitesimal real Ressayre's pair $(\gamma,\Ccal)$ of $\Zcal$.
\item $(\xi,\gamma)\geq (\Psi(\Ccal),\gamma)$ holds for any regular balanced Ressayre's pair $(\gamma,\Ccal)$ of $\Zcal$.
\item $(\xi,\gamma)\geq (\Psi(\Ccal),\gamma)$ holds for any $(\gamma,\Ccal)$ such that 
\begin{enumerate}
 \item[a)] $\gamma$ is rational,
 \item[b)] $\dim_\pgot(\Ccal)-\dim_\pgot(\Zcal)\in\{0,1\}$, 
\item[c)] $({\rm C}_1)$, $({\rm C}_2)$ and $({\rm C}_3)$ hold.
\end{enumerate}
\end{itemize}
The result is still true if we replace condition $c)$ by 
$$
c')\qquad ({\rm C}_1), ({\rm C}_2)\ {\rm and}\ ({\rm C}'_3) \ {\rm hold}.
$$
\end{theorem}

\section{Proof of the main theorem}

We come back to the setting of \S \ref{sec:kahler-hamiltonian}. We will explain how Theorem \ref{theo:main} is a consequence of 
Theorem \ref{theo-facet-kirwan-reel} applied to the K\"ahler-Hamiltonian manifold $M=\T\tU$ equipped with the action of the reductive group 
$L_\C=\tU_\C\times U_\C$.

The real reductive group $G\subset \tG$ are the connected components of the subgroups $U_\C^\sigma\subset \tU_\C^\sigma$. We have a Cartan decomposition $\ggot=\kgot\oplus\pgot$, where $\kgot=\ugot^\sigma$ and $\pgot=i\ugot^{-\sigma}$ (idem for $\tggot$). Let $\agot\subset\tagot$ be maximal abelian subspaces of $\pgot\subset\tpgot$.

\subsection{Admissible elements} 

We see $\tG$ as a submanifold of $\tU_\C$ equipped with the following $\tG\times G$-action : $(\tg,g)\cdot \tilde{x}=\tg\tilde{x} g^{-1}$. 
Let $(\tX,X)$ be an element of  $\tagot\times \agot$. In the following lemma, we describe the manifold $\tG^{(\tilde X,X)}$ of zeroes of the vector 
field $(\tX,X)$ on $\tG$. We denote by $\tG^X\subset \tG$ the subgroup that fixes $X$. 

\begin{lemma}\label{lem:point-fixe}
\begin{itemize}
\item If $\tX\notin \tK X$ , then $\tG^{(\tilde X,X)}=\emptyset$

\item If $\tilde X=\tilde{k} X$ with $\tilde{k}\in \tK$, then $\tG^{(\tilde X,X)}=\tilde{k}\cdot \tG^X$.
\end{itemize}
\end{lemma}

{\em Proof :} The set $\tG^{(\tilde X,X)}$ is non-empty if and only if $\tX$ belongs to the adjoint orbit $\tG X$. The Lemma follows then from the fact that 
the intersection $\tG X\cap \tpgot$ is equal to $\tK X$. $\Box$

\medskip

Let $W_{\tagot}=N_{\tK}(\tagot)/Z_{\tK}(\tagot)$ be the restricted Weyl group. Since both elements $\tX$ and $X$ belongs to $\tagot$,  we notice that if 
$\tX=\tilde{k}X$ for some $\tilde{k}\in \tK$, then there exists $\tilde{w}\in W_{\tagot}$ such that $\tX=\tilde{w}X$. So, the admissible elements relative to the action of $\tG\times G$ on $\tG$ are of the form $(\tilde{w}\gamma,\gamma)$ where $\gamma\in \agot$ is a rational element and  $\tilde{w}\in W_{\tagot}$.

We let $\Sigma(\tggot/\ggot)\subset \agot^*$ denote the set of weights relative to the 
$\agot$-action on $\tggot/\ggot$. If $\gamma\in\agot$, we denote by $\Sigma(\tggot/\ggot)\cap \gamma^\perp$ the subset of weights vanishing against $\gamma$. The main purpose of this section is the proof of the following result.

\begin{proposition}\label{lem:admissible-G-tilde}
Let $G\subset \tG$ be two linear real reductive groups satisfying Assumption \ref{hypothese}. An element $(\tilde{w}\gamma,\gamma)$ is {\em admissible} relatively to the $\tG\times G$-action on $\tG$ if and only if $\gamma$ is {\em rational}  and if 
${\rm Vect}\big(\Sigma(\tggot/\ggot)\cap \gamma^\perp\big)={\rm Vect}\big(\Sigma(\tggot/\ggot)\big)\cap \gamma^\perp$.
\end{proposition}

\medskip 

The proof of Proposition \ref{lem:admissible-G-tilde} is a consequence of the next Lemma.

\begin{lemma}Let $\gamma\in\agot$ and let $\tX_o\in\tagot$ be a regular element.
\begin{enumerate}
\item We have 
$\dim_{\tpgot\times\pgot}\tG=\dim_{\pgot}\tK\tX_o$ 
\ and \ 
$\dim_{\tpgot\times\pgot}\tG^{(\tilde{w}\gamma,\gamma)}=\dim_{\pgot^\gamma}\tK^\gamma\tX_o$.
\item Under Assumption \ref{hypothese}, we have 
\begin{itemize}
\item $\dim_{\tpgot\times\pgot}\tG=\dim Z(\tggot)\cap\pgot$.
\item $\dim_{\pgot^\gamma}\tK^\gamma\tX_o-\dim_{\pgot}\tK\tX_o\in\{0,1\}$ if and only if 
${\rm Vect}\big(\Sigma(\tggot/\ggot)\cap \gamma^\perp\big)$ is equal to ${\rm Vect}\big(\Sigma(\tggot/\ggot)\big)\cap \gamma^\perp$.
\end{itemize}
\end{enumerate}
\end{lemma}

{\em Proof :} For any $x\in \tG$, the stabilizer $(\tpgot\times\pgot)^x$ admits a canonical identification with ${\rm Ad}(x)(\pgot)\cap\tpgot$. If we write 
$x=\tk'e^{\tY}$, with $\tk'\in\tK$ and $\tY\in\tpgot$, we see that ${\rm Ad}(x)(\pgot)\cap\tpgot\simeq{\rm Ad}(e^{\tY})(\pgot)\cap\tpgot=\pgot\cap\tpgot^{\tY}$.
Thus $\dim_{\tpgot\times\pgot}\tG=\min_{\tY\in\tpgot}\dim(\pgot\cap\tpgot^{\tY})$. If we write $\tY=\tk\cdot\tX$, with $\tX\in\tagot$, we have 
${\rm Ad}(\tk)(\tagot)\subset \tk\cdot\tpgot^{\tX}=\tpgot^{\tY}$, and the equality ${\rm Ad}(\tk)(\tagot)=\tpgot^{\tY}$ holds when $\tY$ is regular. So 
$\dim_{\tpgot\times\pgot}\tG=\min_{\tk\in\tK}\dim(\pgot\cap{\rm Ad}(\tk)(\tagot))$. On the other hand, when $\tX_o\in\tagot$ is a regular element it is immediate that 
$\dim_{\pgot}\tK\tX_o=\min_{\tk\in\tK}\dim(\pgot\cap{\rm Ad}(\tk)(\tagot))$. We obtain finally that $\dim_{\tpgot\times\pgot}\tG=\dim_{\pgot}\tK\tX_o$. 
Since $\dim_{\tpgot\times\pgot}\tG^{(\tilde{w}\gamma,\gamma)}= \dim_{\tpgot\times\pgot}\tG^{\gamma}$, the indentity $\dim_{\tpgot\times\pgot}\tG^{(\tilde{w}\gamma,\gamma)}=\dim_{\pgot^\gamma}\tK^\gamma\tX_o$ admits the same proof than before. First one sees that
$$
\dim_{\tpgot\times\pgot}\tG^{\gamma}=\min_{\tk\in\tK^\gamma}\dim(\pgot\cap{\rm Ad}(\tk)(\tagot))=\min_{\tk\in\tK^\gamma}\dim(\pgot^\gamma\cap{\rm Ad}(\tk)(\tagot))
$$
and then we use that $\min_{\tk\in\tK^\gamma}\dim(\pgot^\gamma\cap{\rm Ad}(\tk))=\dim_{\pgot^\gamma}\tK^\gamma\tX_o$.
The first point is settled.

$\tK\tX_o$ is the real part of the coadjoint orbit $\tU\tX_o$. Thanks to Lemma \ref{lem:stabilisateur-min}, we know that exists a subspace $\hgot\subset \pgot$ such that  $\forall x\in \tK\tX_o, \exists k\in K$ such that ${\rm Ad}(k)(\hgot)\subset \pgot^x$, and $\dim(\hgot)=\dim(\pgot^x)$ on a dense open subset of $\tK\tX_o$ : $\hgot$ is the generic $\pgot$-stabilizer on $\tK\tX_o$ and $\dim_{\pgot}(\tK\tX_o)=\dim(\hgot)$. 

Let us check that under Assumption \ref{hypothese},  we have $\hgot=Z(\tggot)\cap\pgot$.
If $\beta'\in\hgot$, then $K(\tK\tX_o)^{\beta'}=\tK\tX_o$. Let us write $\beta'={\rm Ad}(k')(\beta)$ with $\beta\in\agot$ and $k'\in K$. Since $(\tK\tX_o)^{\beta}=\cup_{\tw\in W_{\tagot}}\tK^\beta\tw\tX_o$, the relation $K(\tK\tX_o)^{\beta'}=\tK\tX_o$ implies that a subset $K\tK^\beta\tw\tX_o\subset \tK\tX_o$ has a non empty interior. Let $g=k\tk\tw\in K\tK^\beta\tw$ such that $g\tX_o$ belongs to the interior of $K\tK^\beta\tw\tX_o$: this implies the identity 
\begin{equation}\label{eq:submersion}
\kgot+{\rm Ad}(k)(\tkgot^\beta)+{\rm Ad}(g)(\tkgot^{\tX_o})=\tkgot.
\end{equation}
But ${\rm Ad}(g)(\tkgot^{\tX_o})={\rm Ad}(k\tk)(\tkgot^{\tX_o})\subset{\rm Ad}(k)(\tkgot^{\beta})$, hence (\ref{eq:submersion}) is equivalent to 
$\kgot+\tkgot^\beta=\tkgot$. Now standard computation shows that 
$$
\kgot+\tkgot^\beta=\tkgot\Longleftrightarrow [\beta,\tkgot]\subset\pgot\Longleftrightarrow [\beta,\tpgot]\subset\kgot\Longleftrightarrow [\beta,\tggot]\subset\ggot
\Longleftrightarrow [\beta',\tggot]\subset\ggot.
$$
We have proved that any element $\beta'\in\hgot\subset\pgot$ belongs to the ideal $\mathfrak{b}:=\{X\in\ggot, [X,\tggot]\subset \ggot\}$. Thanks to Assumption \ref{hypothese}, we know that $\mathfrak{b}=Z(\tggot)\cap\ggot$. We have proved that $\hgot\subset Z(\tggot)\cap\pgot$.
On the other hand, it is immediate that $Z(\tggot)\cap\pgot\subset \pgot^x$ for any $x\in\tK\tX_o$. We have proved that the equality 
$Z(\tggot)\cap\pgot=\pgot^x$ holds on a dense open subset of $\tK\tX_o$. The first part of the second point is settled.

Let us concentrated on the last part of the second point. 

If $\gamma\in Z(\tggot)\cap\pgot$ then $\dim_{\pgot^\gamma}\tK^\gamma\tX_o=
\dim_{\pgot}\tK\tX_o$ and ${\rm Vect}\big(\Sigma(\tggot/\ggot)\cap \gamma^\perp\big)={\rm Vect}\big(\Sigma(\tggot/\ggot)\big)\cap \gamma^\perp={\rm Vect}\big(\Sigma(\tggot/\ggot)\big)$.

Suppose now that $\gamma\notin Z(\tggot)\cap\pgot$. We see that $Z(\tggot)\cap\pgot\oplus \R\gamma\subset (\pgot^\gamma)^x$ for any 
$x\in \tK^\gamma\tX_o$, hence
$\dim_{\pgot^\gamma}\tK^\gamma\tX_o-\dim_{\pgot}\tK\tX_o=1$ if and only if $Z(\tggot)\cap\pgot\oplus \R\gamma$ is the generic $\pgot^\gamma$-stabilizer on 
$\tK^\gamma\tX_o$.

Suppose that ${\rm Vect}\big(\Sigma(\tggot/\ggot)\cap \gamma^\perp\big)={\rm Vect}\big(\Sigma(\tggot/\ggot)\big)\cap \gamma^\perp$. For any $\beta\in\agot$, 
we have the equivalences 
\begin{eqnarray*}
[\beta,\tggot^\gamma]\subset \ggot^\gamma\Longleftrightarrow \langle\alpha,\beta\rangle=0,\forall \alpha\in \Sigma(\tggot/\ggot)\cap \gamma^\perp
&\Longleftrightarrow&\beta\in {\rm Vect}\big(\Sigma(\tggot/\ggot)\cap \gamma^\perp\big)^\perp\\
&\Longleftrightarrow&\beta\in \left({\rm Vect}\big(\Sigma(\tggot/\ggot)\big)\cap \gamma^\perp\right)^\perp\\
&\Longleftrightarrow&\beta\in {\rm Vect}\big(\Sigma(\tggot/\ggot))^\perp+ \R\gamma\\
&\Longleftrightarrow&\beta\in Z(\tggot)\cap\pgot\oplus\R\gamma.
\end{eqnarray*}
Let $\hgot_\gamma\subset\pgot^\gamma$ be the generic $\pgot^\gamma$-stabilizer on $\tK^\gamma\tX_o$. Let $\beta'\in \hgot_\gamma$ that we write 
$\beta'=k\beta$ with $\beta\in\agot$ and $k\in K^\gamma$. Since $K^\gamma(\tK^\gamma\tX_o)^{\beta'}=K^\gamma(\tK^\gamma\tX_o)^{\beta}$ is equal to 
$\tK^\gamma\tX_o$, the same arguments than those used before show that $[\beta,\tggot^\gamma]\subset \ggot^\gamma$ and then 
$\beta'=\beta\in Z(\tggot)\cap\pgot\oplus\R\gamma$. We have proved that $\hgot_\gamma=Z(\tggot)\cap\pgot\oplus \R\gamma$ is the generic $\pgot^\gamma$-stabilizer on 
$\tK^\gamma\tX_o$.

Suppose now that ${\rm Vect}\big(\Sigma(\tggot/\ggot)\cap \gamma^\perp\big)\neq{\rm Vect}\big(\Sigma(\tggot/\ggot)\big)\cap \gamma^\perp$. It means that there exists $\beta\notin Z(\tggot)\cap\pgot\oplus \R\gamma={\rm Vect}\big(\Sigma(\tggot/\ggot))^\perp+ \R\gamma$ such that 
$\langle\alpha,\beta\rangle=0,\forall \alpha\in \Sigma(\tggot/\ggot)\cap \gamma^\perp$, i.e. $[\beta,\tggot^\gamma]\subset \ggot^\gamma$. Let consider the subspace 
$\hgot_o=Z(\tggot)\cap\pgot\oplus \R\gamma\oplus\R\beta$. As we have shown before, the inclusion $[\beta,\tggot^\gamma]\subset \ggot^\gamma$ shows that 
the set $K^\gamma(\tK^\gamma\tX_o)^{\hgot_o}=K^\gamma(\tK^\gamma\tX_o)^{\beta}$ has a non-empty interior. Thanks to Lemma \ref{coro:stabilisateur-min}, we can conclude that $\dim_{\pgot^\gamma}\tK^\gamma\tX_o\geq \dim(\hgot_o)$ and then $\dim_{\pgot^\gamma}\tK^\gamma\tX_o-\dim_{\pgot}\tK\tX_o\geq 2$. 
The last part of the second point is proved. $\Box$

\subsection{Balanced Ressayre's pair : infinitesimal conditions} 
Let $(\tilde{w},\gamma)\in W_{\tagot}\times \agot$ such that $\gamma_{\tilde{w}}:=(\tilde{w}\gamma,\gamma)$ is an admissible element relative to the action of 
$\tG\times G$ on $\tG\subset \tU_\C$. 

The complex submanifold of $\tU_\C$ fixed by $\gamma_{\tilde{w}}$ is $C_{\tilde{w},\gamma}:=\tilde{w}\tU^\gamma_\C$ and the real part contained in $\tG$ is 
$\Ccal_{\tilde{w},\gamma}:=C_{\tilde{w},\gamma}\cap \tG=\tilde{w}\tG^\gamma$. The Bialynicki-Birula's complex submanifold associated to 
$C_{\tilde{w},\gamma}$ is $C_{\tilde{w},\gamma}^-:= \tilde{w}\tPbb_\gamma$ where $\tPbb_\gamma$ is the complex parabolic subgroup of 
$\tU_\C$ associated to $\gamma$ (see (\ref{eq:P-gamma-reel})). The parabolic subgroup of 
$\tU_\C\times U_\C$ associated to the weight $\gamma_{\tilde{w}}$ is 
$$
\Pbb_{\tilde{w},\gamma}:=\mathrm{Ad}(\tilde{w})(\tPbb_\gamma)\times \Pbb_\gamma.
$$

Let $\Pbb\subset U_\C$ (resp. $\tPbb\subset \tU_\C$) be the parabolic subgroup associated to the choice of the Weyl chamber $\agot_{+}$
 (resp. $\tagot_{+}$). We can now consider the holomorphic map
$$
{\rm q}_{\tilde{w},\gamma}^\C: (\tPbb\times \Pbb)\times_{(\tPbb\times \Pbb)\cap \Pbb_{\tilde{w},\gamma}} \tilde{w}\tPbb_\gamma\longrightarrow  \tU_\C
$$
that sends $[\tilde{p},p;\tilde{w}\tilde{x}]$ to $\tilde{p}\tilde{w}\tilde{x}p^{-1}$. 
Let us recall the definition of regular balanced Ressayre's pair in this setting.
\begin{definition}
The data $(\gamma_{\tilde{w}},\Ccal_{\tilde{w},\gamma})$ is a regular balanced Ressayre's pair of $\tG$ if there exists $\tau$-invariant 
Zariski open subsets $\Vcal^\C\subset  \tU_\C$ and $\Ucal^\C\subset \tilde{w}\tPbb_\gamma$ such that, $\Vcal^\C$ is  $\tPbb\times \Pbb$-invariant, 
$\Ucal^\C$ is  $(\tPbb\times \Pbb)\cap \Pbb_{\tilde{w},\gamma}$-invariant and intersects $\tilde{w}\tU^\gamma_\C$, 
annd the map ${\rm q}_{\tilde{w},\gamma}^\C$ defines a diffeomorphism 
\begin{equation}\label{eq:diffeo-balanced}
(\tPbb\times \Pbb)\times_{(\tPbb\times \Pbb)\cap \Pbb_{\tilde{w},\gamma}}\Ucal^\C\simeq \Vcal^\C.
\end{equation}
\end{definition}

As we have explain in \S \ref{sec:balanced-RP}, the set $\Ucal^\C\cap \tilde{w}\tG^\gamma$ is nonempty and the  diffeomorphism (\ref{eq:diffeo-balanced}) insures that the infinitesimal map $\rho^{\gamma_{\tilde{w}}}_z: \tngot^{\tilde{w}\gamma>0}\times \ngot^{\gamma>0}\longrightarrow (\T_z \tG)^{\gamma_{\tw}>0}$ is an isomorphism for any $z\in \Ucal^\C\cap\tilde{w}\tG^\gamma$.

\begin{lemma}\label{lem:A-B-balanced}
Let $(\gamma_{\tilde{w}},\Ccal_{\tilde{w},\gamma})$ be a regular balanced Ressayre's pair of $\tG$. Then we have 
$$
({\rm C}_1)\qquad \dim(\tngot^{\tilde{w}\gamma>0})+\dim (\ngot^{\gamma>0})= \dim (\tggot^{\gamma>0}),
$$
$$
({\rm C}_2)\qquad {\rm Tr}_{\gamma}(\ngot^{\gamma>0})+{\rm Tr}_{\tw\gamma}(\tngot^{\tw\gamma>0})={\rm Tr}_{\gamma}(\tggot^{\gamma>0}).
$$
\end{lemma}

{\em Proof :} Let $z\in \tilde{w}\tG^\gamma$ such that $\rho^{\gamma_{\tilde{w}}}_z: \tngot^{\tilde{w}\gamma>0}\times \ngot^{\gamma>0}\longrightarrow (\T_z \tG)^{\gamma_{\tw}>0}$ is an isomorphism. Since $\rho^{\gamma_{\tilde{w}}}_z$ commutes with the infinitesimal action of $\gamma_{\tilde{w}}$, it implies that 
$\dim(\tngot^{\tilde{w}\gamma>0})+\dim (\ngot^{\gamma>0})= \dim ((\T_z \tG)^{\gamma_{\tw}>0})$ and 
${\rm Tr}_{\gamma_{\tilde{w}}}\left(\tngot^{\tilde{w}\gamma>0}\times \ngot^{\gamma>0}\right)={\rm Tr}_{\gamma_{\tilde{w}}}((\T_z \tG)^{\gamma_{\tw}>0})$.

Since $(\T_z \tG)^{\gamma_{\tw}>0}\simeq \tggot^{\gamma>0}$, relations $({\rm C}_1)$ and $({\rm C}_2)$ follows.  $\Box$

\subsection{Cohomological conditions} 

To any $\tw\in W_{\tagot}$, we associate the Schubert cell $\tilde{\Xgot}^{o}_{\tilde{w},\gamma}=\tPbb[\tilde{w}]=\tBbb[\tilde{w}]\subset \tFcal_\gamma$ 
and the Schubert variety $\tilde{\Xgot}_{\tilde{w},\gamma}=\overline{\tilde{\Xgot}^{o}_{\tilde{w},\gamma}}$ (see \S \ref{sec:bruhat}). 
Let us denote $\Xgot^o_\gamma:=\B[e]\subset \Fcal_\gamma$ and the corresponding Schubert variety 
$\Xgot_\gamma=\overline{\Xgot^o_\gamma}$.

Let $[\tilde{\Xgot}_{\tilde{w},\gamma}]\in H^{2n_{\tilde{w},\gamma}}(\tFcal_\gamma,\Z)$ and $[\Xgot_{\gamma}]\in H^{2n_{\gamma}}(\Fcal_\gamma,\Z)$ be the cycle classes defined by these algebraic subvarieties. Here $n_{\tilde{w},\gamma}={\rm codim}_\C(\tilde{\Xgot}_{\tilde{w},\gamma})$ and 
$n_{\gamma}={\rm codim}_\C(\Xgot_{\gamma})$.

A standard computation shows that condition $({\rm C}_1)$ of Lemma \ref{lem:A-B-balanced}  is equivalent to 
${\rm codim}_\C(\tilde{\Xgot}_{\tilde{w},\gamma})+{\rm codim}_\C(\Xgot_{\gamma})=\dim_\C(\Fcal_\gamma)$.

If $(\gamma_{\tilde{w}},\Ccal_{\tilde{w},\gamma})$ is a  balanced Ressayre's pair of $\tG$, then the map ${\rm q}_{\tilde{w},\gamma}^\C$ 
satisfies condition $({\rm C}_3)$ : $\{\tg\in\tU_\C, \sharp \{({\rm q}_{\tilde{w},\gamma}^\C)^{-1}(\tg^{-1})\}=1\} \ {\rm contains\ a\ Zariski\ open\ subset\ of}\ \tU_\C$.

Let us consider the map
$$
\pi_{\tilde{w},\gamma}: (\tPbb\times \Pbb)\times_{(\tPbb\times \Pbb)\cap \Pbb_{\tilde{w},\gamma}} \tw\tPbb_\gamma
\longrightarrow \tFcal_\gamma\times \Fcal_\gamma
$$
 that sends $[\tp, p; \tw x]$ to $([\tp\tw],[p])$.

The following facts are standard (see Proposition 4.3 in \cite{Brion-Bourbaki} or Section 6.1 in \cite{pep-real-RP1}). For any $\tg\in\tU_\C$,  the map $\pi_{\tilde{w},\gamma}$ is injective on the fiver $({\rm q}_{\tilde{w},\gamma}^\C)^{-1}(\tg^{-1})$, and the set \break $\pi_{\tilde{w},\gamma}(({\rm q}_{\tilde{w},\gamma}^\C)^{-1}(\tg^{-1}))$ is equal to 
 $$
  \Xgot^o_\gamma\cap \iota^{-1}(\tg\tilde{\Xgot}_{\tw,\gamma}^o)\simeq \{(x,y)\in \tPbb[\tw]\times \Pbb[e], \tg x= \iota(y)\}\subset \tFcal_\gamma\times \Fcal_\gamma.
 $$
 Here $\iota :\Fcal_\gamma\to\tFcal_\gamma$ is the canonical embedding. Thus conditions $({\rm C}_3)$ and $({\rm C}'_3)$ becomes
 $$
({\rm C}_3) \qquad \left\{\tg\in\tU_\C, \sharp \{\Xgot^o_\gamma\cap \iota^{-1}(\tg\tilde{\Xgot}_{\tw,\gamma}^o)\}=1\right\} \ {\rm contains\ a\ Zariski\ open\ subset\ of}\ \tU_\C.
$$
$$
({\rm C}'_3)\qquad \left\{\tg\in\tU_\C, \Xgot^o_\gamma\cap \iota^{-1}(\tg\tilde{\Xgot}_{\tw,\gamma}^o)\neq\emptyset\right\} \ {\rm contains\ a\ Zariski\ open\ subset\ of}\ \tU_\C.
$$

According to a transversality theorem due to Kleiman, the intersection $\Xgot_\gamma\cap \iota^{-1}(\tg\tilde{\Xgot}_{\tw,\gamma})$ is transverse for any $\tg$ in a nonempty Zariski open subset of $\tU_\C$; moreover, $\Xgot^o_\gamma\cap \iota^{-1}(\tg\tilde{\Xgot}_{\tw,\gamma}^o)$ is dense in this intersection. We deduce that 
\begin{enumerate}
\item $({\rm C}_3)$ is equivalent to asking that the cohomology class of $\Xgot_\gamma\cap \iota^{-1}(\tg\tilde{\Xgot}_{\tw,\gamma})$ is equal to the cohomology class of a point : $[\Xgot_{\gamma}]\cdot \iota^*([\tilde{\Xgot}_{\tilde{w},\gamma}])= [pt]$ in $H^*(\Fcal_\gamma,\Z)$.
\item $({\rm C}_1)$ and $({\rm C}'_3)$ are equivalent to asking that 
$[\Xgot_{\gamma}]\cdot \iota^*([\tilde{\Xgot}_{\tilde{w},\gamma}])= \ell[pt]$ with $\ell\geq 1$ in $H^*(\Fcal_\gamma,\Z)$.
\end{enumerate}
Let us notice that $[\Xgot_{\gamma}]\cdot \iota^*([\tilde{\Xgot}_{\tilde{w},\gamma}])= [pt]$ implies that  ${\rm codim}_\C(\tilde{\Xgot}_{\tilde{w},\gamma})+{\rm codim}_\C(\Xgot_{\gamma})=\dim_\C(\Fcal_\gamma)$, hence $({\rm C}_1)$ is a consequence of $({\rm C}_3)$.

\subsection{Conclusion} 
Let us consider the gradient map $\Psi : \tG\simeq \tK\times\tpgot \to \tpgot\times\pgot, (\tk,\tX) \mapsto (-\tk\tX,\pi(\tX))$.

Let $(\tilde{w},\gamma)\in W_{\tagot}\times \agot$ and $\Ccal_{\tilde{w},\gamma}=\tilde{w}\tG^\gamma$. In the previous sections, we have proved the following
\begin{itemize}
 \item $\gamma_{\tilde{w}}:=(\tilde{w}\gamma,\gamma)$ is admissible if and only if ${\rm Vect}\big(\Sigma(\tggot/\ggot)\cap \gamma^\perp\big)={\rm Vect}\big(\Sigma(\tggot/\ggot)\big)\cap \gamma^\perp$.
 \item Condition $({\rm C}_2)$ for $(\gamma_{\tilde{w}},\Ccal_{\tilde{w},\gamma})$ is equivalent to 
 $ {\rm Tr}_{\gamma}(\ngot^{\gamma>0})+{\rm Tr}_{\tw\gamma}(\tngot^{\tw\gamma>0})={\rm Tr}_{\gamma}(\tggot^{\gamma>0})$.
\item Condition $({\rm C}_3)$ for $(\gamma_{\tilde{w}},\Ccal_{\tilde{w},\gamma})$ is equivalent to $[\Xgot_{\gamma}]\cdot \iota^*([\tilde{\Xgot}_{\tilde{w},\gamma}])= [pt]$  in $H^*(\Fcal_\gamma,\Z)$.
\item $({\rm C}_1)$ is a consequence of $({\rm C}_3)$.
\item Conditions $({\rm C}_1)$ and $({\rm C}'_3)$ for $(\gamma_{\tilde{w}},\Ccal_{\tilde{w},\gamma})$ are equivalent to 
$[\Xgot_{\gamma}]\cdot \iota^*([\tilde{\Xgot}_{\tilde{w},\gamma}])= \ell[pt]$ with $\ell\geq 1$ in $H^*(\Fcal_\gamma,\Z)$.
\end{itemize}

We notice that $(\Psi(\Ccal_{\tilde{w},\gamma}),\gamma_{\tilde{w}})=0$. Thus Theorem \ref{theo-facet-kirwan-reel} gives us the following
\begin{theorem}\label{theo:main-B}Let $(\tilde{\xi},\xi)\in\tagot_+\times \agot_+$. 
We have  $-K\xi\subset \pi\big(\tK\txi\big)$ if and only if 
\begin{equation}\label{eq:inequality-polytope-conclusion}
(\tilde{\xi},\tilde{w}\gamma)+(\xi,\gamma)\geq 0
\end{equation}
for any $(\gamma,\tilde{w})\in\agot\times W_{\tagot}$ satisfying the following properties:
\begin{enumerate}
\item[a)] $\gamma$ is rational and ${\rm Vect}\big(\Sigma(\tggot/\ggot)\cap \gamma^\perp\big)={\rm Vect}\big(\Sigma(\tggot/\ggot)\big)\cap \gamma^\perp$.
\item[b)] $[\Xgot_{\gamma}]\cdot \iota^*([\tilde{\Xgot}_{\tilde{w},\gamma}])= [pt]$ in $H^*(\Fcal_\gamma,\Z)$.
\item[c)] ${\rm Tr}_{\gamma}(\ngot^{\gamma>0})+{\rm Tr}_{\tw\gamma}(\tngot^{\tw\gamma>0})={\rm Tr}_{\gamma}(\tggot^{\gamma>0})$.
\end{enumerate}
The result still holds if we replace condition $b)$ with 
$$
b')\qquad [\Xgot_{\gamma}]\cdot \iota^*([\tilde{\Xgot}_{\tilde{w},\gamma}])= \ell[pt]\ {\rm with}\ \ell\geq 1\ {\rm in}\ H^*(\Fcal_\gamma,\Z).
$$
\end{theorem}

\medskip

We have to make some small modifications in the previous theorem in order to obtain Theorem \ref{theo:main}. The longest element $w_0\in W$, which is invariant by $\sigma$, can be seen
as the unique element of $W_\agot$ such that $w_0(\agot_+)=-\agot_+$. Any elements $\xi\in \agot_+$ and $\gamma \in\agot$ can be written 
$\xi=-w_0(\xi')$ and $\gamma=w\gamma'$ where $\xi'\in \agot_+$, $\gamma'\in -\agot_+$ and $w\in W_\agot$. Since $W_\agot\gamma=K\gamma\cap\agot$ is contained in 
$W_{\tagot}\gamma=\tK\gamma\cap\tagot$, for any $(w,\tw)\in W_\agot\times W_{\tagot}$, there exists $\tw'\in W_{\tagot}$ such that 
$\tw w\gamma'=\tw'\gamma'$.

We leaves the reader verify the following
\begin{itemize}
\item $-K\xi\subset \pi\big(\tK\txi\big)$ is equivalent to $K\xi'\subset \pi\big(\tK\txi\big)$
\item (\ref{eq:inequality-polytope-conclusion}) is equivalent to $(\tilde{\xi},\tilde{w}'\gamma')\geq (\xi',w_0w\gamma')$
\item ${\rm Vect}\big(\Sigma(\tggot/\ggot)\cap \gamma^\perp\big)={\rm Vect}\big(\Sigma(\tggot/\ggot)\big)\cap \gamma^\perp$ is equivalent to 
${\rm Vect}\big(\Sigma(\tggot/\ggot)\cap (\gamma')^\perp\big)={\rm Vect}\big(\Sigma(\tggot/\ggot)\big)\cap (\gamma')^\perp$.
\item Condition $b)$ is equivalent to 
$[\Xgot_{w,\gamma'}]\cdot \iota^*([\tilde{\Xgot}_{\tilde{w}',\gamma'}])= [pt]$ in $H^*(\Fcal_{\gamma'},\Z)$.
\item Condition $c)$ is equivalent to 
${\rm Tr}_{w\gamma'}(\ngot^{w\gamma'>0})+{\rm Tr}_{\tw\gamma'}(\tngot^{\tw\gamma'>0})={\rm Tr}_{\gamma'}(\tggot^{\gamma'>0})$.
\end{itemize}

We see then that Theorem \ref{theo:main} is a reformulation of Theorem \ref{theo:main-B}.

\section{Singular value inequalities for matrix sums}

Let $p\geq q\geq 1$ and $n=p+q$. We consider the canonical action of the group $K:=U(p)\times U(q)$ on the vector space 
$M_{p,q}(\C)$ of complex matrices of size $p\times q$. Let $\agot\subset M_{p,q}(\C)$ be the real subspace consisting of all matrices of the form 
$$
Y(x)=\begin{pmatrix}
0&\cdots & x_1\\
\vdots &\reflectbox{$\ddots$}& \vdots\\
x_q&\cdots & 0\\
0&\cdots & 0\\
\vdots &\ddots& \vdots\\
0&\cdots & 0\\
\end{pmatrix},\, x\in\R^q.
$$ 
The chamber $\agot_+:=\{Y(x), x_1\geq x_2\geq \cdots\geq x_q\geq 0\}$ parametrizes the $K$-orbits on $M_{p,q}(\C)$~: in other words a matrix $X$ belongs to the orbit 
$K\cdot Y(x)$ if and only $x$ is equal to the singular spectrum of $X$ (see \S \ref{sec:singular-horn}). Thus we can consider 
$$
\sing(p,q):=\left\{(a,b,c)\in (\R^q_{++})^3, K\cdot Y(c)\subset K\cdot Y(a)+K\cdot Y(b)\right\},
$$
where $\R^q_{++}=\{x=(x_1\geq x_2\geq \cdots\geq x_q\geq 0)\}$. Notice that $K\cdot Y(x)=-K\cdot Y(x)$ for any $x\in \R^q_{++}$, hence $\sing(p,q)$ can be 
defined with the symmetric condition $0\in K\cdot Y(a)+K\cdot Y(b)+K\cdot Y(c)$.

We need to recall some notations introduced in \S \ref{sec:singular-horn}. For any $r\in [q]$, let 
$\Pcal^{p,q}_r$ be the collection of subsets $I\subset [n]$ of cardinal $r$, satisfying $I\cap I^o=\emptyset$ and $I\cap \{q+1,\ldots, p\}=\emptyset$, where $I^o=\{n+1-k, k\in I\}$. 
We see that any $I\in \Pcal^{p,q}_r$ is equal to the disjoint union $I_+\coprod I_{-}^o$ where $I_{+},I_{-}$ are disjoint subsets of $[q]$. A subset $A\subset [n]$ is 
polarized if it admits a decomposition $A=A_+ \coprod A_-$ into disjoints subsets. If  $A_+$ and $A_-$ are both of cardinal $r$, we say that $A$ is balanced. We denote by $\Pol^n_{2r}$ the collection of balanced polarized subsets of $[n]$ of cardinal $2r$.

To any subset $I\in\Pcal^{p,q}_r$, we associate the balanced polarized subset $\widehat{I}=I\coprod I^o\in  \Pol^n_{2r}$ and the cohomology class 
$[\Xgot_{\widehat{I}}]\in H^*(\F(r,n-r;n),\Z)$. The proof of the following result is given is \S \ref{sec:preuve-singular-horn}.

\begin{theorem}\label{theo:singular-p-q-section}
An element $(x,y,z)\in(\R^q_{++})^3$ belongs to $\sing(p,q)$ if and only if 
$$
(\star)_{I_\pm,J_\pm,K_\pm}\hspace{2cm} |\,x\,|_{I_+}+|\,y\,|_{J_+} +|\, z\,|_{K_+}\leq |\,x\,|_{I_{-}}+|\,y\,|_{J_{-}} +|\, z\,|_{K_{-}}
$$
holds for any $r\in [q]$, for any $I,J,K\in\Pcal^{p,q}_r$ such that 
\begin{equation}\label{eq:condition-cohomologique-2-flag}
[\Xgot_{\widehat{I}}]\odot_0[\Xgot_{\widehat{J}}]\odot_0[\Xgot_{\widehat{K}}]= [pt]\quad {\rm in}\quad  H^*(\F(r,n-r;n),\Z).
\end{equation}
The result still holds if we replace (\ref{eq:condition-cohomologique-2-flag}) by the weaker condition 
$[\Xgot_{\widehat{I}}]\odot_0[\Xgot_{\widehat{J}}]\odot_0[\Xgot_{\widehat{K}}]=\ell [pt]$ with $\ell\geq 1$.
\end{theorem}

The following examples are detailed in \S \ref{sec:example-horn-1-2-3}.

\begin{example}
Let us consider the case $p\geq q=1$. Thus $(a,b,c)\in(\R_+)^3$ belongs to $\sing(p,1)$ if and only if the Weyl inequalities holds :
$a+b\geq c$, $a+c\geq b$, and $b+c\geq a$.
\end{example}

\begin{example}
Let us consider the case $p\geq q=2$. Thus $(a,b,c)\in(\R^2_{++})^3$ belongs to $\sing(p,2)$ if and only if the following 18 inequalities holds 
\begin{enumerate}
\item the Weyl inequalities
\begin{itemize}
\item $a_1+ b_1\geq  c_1$ \quad (and 2 permutations),
\item $a_1+ b_2\geq  c_2$ \quad (and 5 permutations),
\end{itemize}
\item the Lidskii inequalities
\begin{itemize}
\item $a_1+a_2+b_1+b_2\geq c_1+c_2$\quad (and 2 permutations),
\end{itemize}
\item the signed Lidskii inequalities
\begin{itemize}
\item $a_1 + a_2 + b_1 - b_2 \geq c_1 - c_2$\quad (and 5 permutations).
\end{itemize}
\end{enumerate}
\end{example}

\subsection{Multiplicative formulas in $H^*(\F(r,n-r;n),\Z)$}\label{sec:multiplicative}

In this section, we use some multiplicative formula obtained by Ressayre \cite{Ressayre11} and Richmond \cite{Richmond12} 
for the structure constants in the cohomology of flag varieties to simplify condition (\ref{eq:condition-cohomologique-2-flag}). We are interested by their multiplicative formula in the case of two steps flag varieties $\F(r,n-r;n)$. The Schubert varieties of $\F(r,n-r;n)$ are parametrized by the set $\Pol^n_{2r}$ of balanced polarized subsets of $[n]$ of cardinal $2r$. If $A=A_+ \coprod A_-\in \Pol^n_{2r}$, the subset $(A_+)^c\subset [n]$ is of cardinal $n-r$ : $(A_+)^c=\{u_1<\cdots<u_{n-r}\}$. Since $A_{-}$ is contained in 
$(A_+)^c$ we can define the following subsets of cardinal $r$ : $A'=\{i\in [n-r], u_i\in A_{-}\}$ and $A''=\{k\in [n-r], n-r+1-k\in A'\}$. Therefore, we can associate the cycles classes $[\Xgot_{A_+}]\in H^*(\G(r,n),\Z)$ and $[\Xgot_{A''}]\in H^*(\G(r,n-r),\Z)$ to any $A\in\Pol^n_{2r}$.

The following result is a particular case of a general multiplicative formula obtained by Ressayre \cite{Ressayre11} and Richmond \cite{Richmond12}.

\begin{proposition}\label{prop:multiplicative}
Let $A,B,C\in\Pol^n_{2r}$. The relation $ [\Xgot_A]\odot_0[\Xgot_B]\odot_0 [\Xgot_C]=\ell[pt]$, with $\ell>0$, holds in $H^*(\F(r,n-r;n),\Z)$ if and only if we have 
\begin{enumerate}
\item $[\Xgot_{A_+}]\cdot[\Xgot_{B_+}]\cdot[\Xgot_{C_+}]= \ell'[pt]$ with $\ell'>0$ in $H^*(\G(r,n),\Z)$,
\item $[\Xgot_{A''}]\cdot[\Xgot_{B''}]\cdot[\Xgot_{C''}]= \ell''[pt]$ with $\ell''>0$ in $H^*(\G(r,n-r),\Z)$,
\item $|A_+|+|B_+|+|C_+|= |A_-|+|B_-|+|C_-| +r(n-r)$.
\end{enumerate}
Moreover we have $\ell=\ell'\ell''$.
\end{proposition}

When we work with $I,J,K\in\Pcal^{p,q}_r$ and the associated balanced polarized subsets $\widehat{I}=I\coprod I^o, \widehat{J}=J\coprod J^o$, and 
$\widehat{K}=K\coprod K^o$, we see that
\begin{eqnarray*}
|\widehat{I}_+| + |\widehat{J}_+|+|\widehat{K}_+|= |\widehat{I}_-| + |\widehat{J}_-|+|\widehat{K}_-| +r(n-r)  & \Longleftrightarrow&\\
 |I| + |J|+|K| =2r(n-r)+\tfrac{3}{2}r(r+1) & \Longleftrightarrow& \\
 {\rm codim}(\Xgot_{I})+{\rm codim}(\Xgot_{J})+{\rm codim}(\Xgot_{K})=r(n-r).
\end{eqnarray*}
Therefore condition {\em 3.} is a consequence of condition {\em 1.} in Proposition \ref{prop:multiplicative} when working with $\widehat{I},\widehat{J},\widehat{K}\in\Pol^n_{2r}$

If $X\in\Pcal^{p,q}_r$, we define $\widetilde{X}\in\Pcal^{n-r}_r$ as equal to $(\widehat{X})''$. Thanks to Proposition \ref{prop:multiplicative}, 
we can state another version of Theorem \ref{theo:singular-p-q}.

\begin{theorem} \label{theo:singular-p-q-version-2}
An element $(x,y,z)\in(\R^q_{++})^3$ belongs to $\sing(p,q)$ if and only if 
$(\star)_{I_\pm,J_\pm,K_\pm}$ holds for any $r\in [q]$, for any $I,J,K\in\Pcal^{p,q}_r$ satisfying both conditions :
\begin{equation}\label{eq:singular-1}
[\Xgot_{I}] \cdot [\Xgot_{J}] \cdot [\Xgot_{K}] = [pt] \quad in\quad  H^*(\G(r,n),\Z),
\end{equation}
\begin{equation}\label{eq:singular-2}
[\Xgot_{\widetilde{I}}] \cdot [\Xgot_{\widetilde{J}}] \cdot [\Xgot_{\widetilde{K}}] =  [pt]\quad in \quad  H^*(\G(r,n-r),\Z).
\end{equation}
The result is still true if we work with the weaker conditions :
$[\Xgot_{I}] \cdot [\Xgot_{J}] \cdot [\Xgot_{K}] =\ell' [pt]$ and $[\Xgot_{\widetilde{I}}] \cdot [\Xgot_{\widetilde{J}}] \cdot [\Xgot_{\widetilde{K}}] = \ell'' [pt]$ with $\ell',\ell''\geq 1$.
\end{theorem}

\subsection{Inequalities of Weyl, Lidskii and Thompson }\label{sec:WLT-inequalities}

The purpose of this section is to  give some examples of triplets $(I,J,K)$ satisfying (\ref{eq:singular-1}) and (\ref{eq:singular-2}) which will allow us to 
recover the classical  inequalities of Weyl, Lidskii and Thompson.

We need to recall a standard feature of Littlewood-Richardson coefficients. Let  $m\geq 2$ and $r\in[m]$. For any subset $X=\{x_1<\cdots<x_r\}\subset [q]$ 
of cardinal $r$, we associate several objects : the partition $\lambda(X)=(\lambda_1\geq\cdots\geq \lambda_r)$ where $\lambda_a=m-r+a-x_a$ ; the polynomial irreducible representation $V_{\lambda(X)}$ of $U(r)$ and the cycle class $[\Xgot_X]\in H^{2n(X)}(\G(r,m-r),\Z)$, where $n(X)=|\lambda(X)|$. A classical result of Lesieur \cite{Lesieur} tell us that for any subsets  $A,B,C\subset[m]$ of cardinal $r$, if $|\lambda(A)|+|\lambda(B)|+|\lambda(C)|=r(m-r)$ we have 
$$
[\Xgot_A]\cdot[\Xgot_B]\cdot[\Xgot_C]=c(A,B,C)[pt]
$$
where $c(A,B,C)=\dim (V_{\lambda(A)}\otimes V_{\lambda(B)}\otimes V_{\lambda(C)}\otimes \det^{r-m})^{U(r)}$. Since $|\lambda(A)|+|\lambda(B)|+|\lambda(C)|=r(m-r)$ 
when $c(A,B,C)\neq 0$, we obtain the following equivalence for any $ \ell\geq 1$ :
$$
[\Xgot_A]\cdot[\Xgot_B]\cdot[\Xgot_C]=\ell[pt]\quad \Longleftrightarrow\quad  
\dim (V_{\lambda(A)}\otimes V_{\lambda(B)}\otimes V_{\lambda(C)}\otimes \det{}^{r-m})^{U(r)}=\ell.
$$

\subsubsection*{Weyl inequalities}

Let $m\geq 2$. To any element $a\in [m]$ we associate the Schubert variety $\Xgot_a=\overline{B_m\cdot (\C e_a)}\subset \G(1,m)$ of codimension $m-a$. Lesieur's result  leads to the following simple criterion.

\begin{lemma}
Let $i,j,k\in [m]$. Then $[\Xgot_i] \cdot [\Xgot_j] \cdot [\Xgot_k]= [pt]$  in $H^*(\G(1,m),\Z)$ if and only if $i+j+k=2m+1$.
\end{lemma}

Let $i,j\in [q]$ such that $i+j\leq q+1$. Consider the following triplet of $\Pcal^{p,q}_1$ : $I=\emptyset \coprod \{n+1-i\}$, 
$J=\emptyset \coprod \{n+1-j\}$ and $K=\{i+j-1\}  \coprod\emptyset$.
On one side, the cycles classes in $H^*(\G(1,n),\Z)$ are $[\Xgot_I]=[\Xgot_{n+1-i}]$, $[\Xgot_J]=[\Xgot_{n+1-j}]$ and $[\Xgot_K]=[\Xgot_{i+j-1}]$. On the other side,  
the cycles classes in $H^*(\G(1,n-1),\Z)$ are $[\Xgot_{\widetilde{I}}]=[\Xgot_{n-i}]$, $[\Xgot_{\widetilde{J}}]=[\Xgot_{n-j}]$ and $[\Xgot_{\widetilde{K}}]=[\Xgot_{i+j-1}]$. Thanks to the previous Lemma, we see that Conditions (\ref{eq:singular-1}) and (\ref{eq:singular-2}) are satisfied. The corresponding inequalities are the {\em Weyl inequalities} 
\cite{Weyl} :
$$
\tau_{i+j-1}(A+B)\leq \tau_i(A)+\tau_j(B),\qquad \forall A,B\in M_{p,q}(\C).
$$
whenever $i+j\leq q+1$.

\subsubsection*{Lidskii  inequalities}

Let us consider the element $X^n_r:=\emptyset\,\coprod \,\{n-r+1,\ldots,n\}\in \Pcal^{p,q}_r$ : here $\widetilde{X}^n_r= \{n-2r+1,\ldots,n-r\}\in \Pcal^{n-r}_r$. 
The Schubert varieties $\Xgot_{X^n_r}$ and $\Xgot_{\widetilde{X}^n_r}$ are respectively equal to the whole manifolds $\G(r,n)$ and $\G(r,n-r)$, hence $[\Xgot_{X^n_r}] = 1$ and 
$[\Xgot_{\widetilde{X}^n_r}] = 1$. Recall that the identity $[\Xgot_{I}] \cdot [\Xgot_{J}] = [pt]$ holds in  $H^*(\G(r,n),\Z)$ if and only if $J=I^o$. 

Thus we obtain the following 
\begin{lemma}
Let $I\in\Pcal^{p,q}_r$. The triplet $(X^n_r,I,I^o)$ satisfies (\ref{eq:singular-1}) and (\ref{eq:singular-2}) if and only if $\widetilde{I^o}=(\widetilde{I}\,)^o$.
\end{lemma}

We leave it to the reader to verify that the identity $\widetilde{I^o}=(\widetilde{I}\,)^o$ is satisfied in the following two cases:
\begin{enumerate}
\item $I=\{i_1<\cdots<i_r\}\coprod \emptyset$. The corresponding inequalities are the {\em Lidskii inequalities} \cite{Lidskii}:
$$
\sum_{k=1}^r\tau_{i_k}(A+B)\leq \sum_{k=1}^r\tau_{i_k}(B) +\sum_{k=1}^r \tau_{k}(A),\qquad \forall A,B\in M_{p,q}(\C).
$$
\item $I=\{i_1<\cdots <\widehat{i_s}< \cdots <i_r\}\coprod \{n+1 -i_s\}$. The corresponding inequalities will be referred to as the {\em signed Lidskii inequalities} : 
$$
\sum_{k\neq s}\tau_{i_k}(A+B)-\tau_{i_s}(A+B)\leq \sum_{k\neq s}\tau_{i_k}(B) -\tau_{i_s}(B)+\sum_{k=1}^r \tau_{k}(A),\qquad \forall A,B\in M_{p,q}(\C).
$$
\end{enumerate}

\subsubsection*{Thompson inequalities}

If $X=\{x_1<\cdots<x_r\}\subset [q]$ is a subset of cardinal $r\geq 1$, we define the partition $\mu_X=(x_r-r\geq \cdots\geq x_1-1)$ and the irreducible representation $V_{\mu_X}$ of $U(r)$ with highest $\mu_X$. Recall that $V_{\mu_X}^*=V_{\mu^*_X}$ where $\mu^*_X=(1-x_1\geq \cdots\geq r-x_r)$.

Let $r\in [q]$. Let us consider three subsets of $[q]$ of cardinal $r$ : $I_-:=\{i_1<\cdots<i_r\}$, $J_-:=\{j_1<\cdots<j_r\}$ and 
$K_+=\{k_1<\cdots<k_r\}$. For any $p\geq q$, we define three elements of $\Pcal^{p,q}_r$ as follows : $I^p=\emptyset \coprod \{p+q+1-i, i\in I_-\}$,  
$J^p=\emptyset \coprod \{p+q+1-j, i\in J_-\}$ and $K^p=K_+\coprod \emptyset$.

\begin{lemma} The triplet $(I^p,J^p,K^p)$ defined above satisfies (\ref{eq:singular-1}) and (\ref{eq:singular-2}) if and only if 
$\dim(V_{\mu_{I_-}}\otimes V_{\mu_{J_-}}\otimes V_{\mu_{K_+}}^*)^{U(r)}=1$. If $(V_{\mu_{I_-}}\otimes V_{\mu_{J_-}}\otimes V_{\mu_{K_+}}^*)^{U(r)}\neq \{0\}$, then 
the weak versions of (\ref{eq:singular-1}) and (\ref{eq:singular-2}) are true.
\end{lemma}

{\em Proof :} A direct computation gives  : $\lambda(I^p)=\lambda(\widetilde{I^p})=\mu_{I_-}$, $\lambda(J^p)=\lambda(\widetilde{J^p})=\mu_{J_-}$ and 
$\lambda(K^p)=(n-r) 1^r +\mu_{K_+}^*$, $\lambda(\widetilde{K^p})=(n-2r) 1^r +\mu_{K_+}^*$. Using Lesieur's result one sees that, for any $\ell\geq 1$, the following statements are equivalent :
\begin{itemize}
\item $[\Xgot_{I^p}] \cdot [\Xgot_{J^p}] \cdot [\Xgot_{K^p}] = \ell[pt]$ in $H^*(\G(r,n),\Z)$,
\item $[\Xgot_{\widetilde{I^p}}] \cdot [\Xgot_{\widetilde{J^p}}] \cdot [\Xgot_{\widetilde{K^p}}] = \ell [pt]$ in $H^*(\G(r,n-r),\Z)$.
\item $\dim(V_{\mu_{I_-}}\otimes V_{\mu_{J_-}}\otimes V_{\mu_{K_+}}^*)^{U(r)}=\ell$.
\end{itemize}
$\Box$

\begin{corollary} We have 
$$
\sum_{k\in K_+}\tau_{k}(A+B)\leq \sum_{j\in J_-}\tau_{j}(B) +\sum_{i\in I_-} \tau_{i}(A),\quad \forall A,B\in M_{p,q}(\C)
$$
 if $(V_{\mu_{I_-}}\otimes V_{\mu_{J_-}}\otimes V_{\mu_{K_+}}^*)^{U(r)}\neq \{0\}$.
\end{corollary}

Let us consider the situation where $k_\alpha = i_\alpha+j_\alpha-\alpha$,  $\forall\alpha\in [r]$. Theses relations imply that 
$\mu_{K_+}=\mu_{I_-}+\mu_{J_-}$, and a standard fact of representation theory tells us that 
$\dim(V_{\mu_{I_-}}\otimes V_{\mu_{J_-}}\otimes V_{\mu_{I_-}+\mu_{J_-}}^*)^{U(r)}=1$. We obtain here Thompson's inequalities \cite{Thompson} 
$$
\sum_{\alpha=1}^r\tau_{i_\alpha+j_\alpha-\alpha}(A+B)\leq \sum_{\alpha=1}^r\tau_{j_\alpha}(B) +\sum_{\alpha=1}^r \tau_{i_\alpha}(A),\quad \forall A,B\in M_{p,q}(\C),
$$
which are valid for any subsets $\{i_1<\cdots<i_r\}$ and $\{j_1<\cdots<j_r\}$ satisfying $i_r + j_r-r\leq q$.

\subsection{Horn conditions}\label{sec:horn-condition}

In the previous sections we parameterized the inequality $(\star)_{I_\pm,J_\pm,K_\pm}$ in terms of subsets $I,J,K$ leaving in $\Pcal^{p,q}_r$. 
To simplify the presentation, we will choose a parameterization using directly the $I_\pm,J_\pm,K_\pm$.

The cardinal of a finite set $Z$ is denoted by $\sharp Z$. Let start from a polarized subset $X_\bullet=X_+\coprod X_-$ of $[q]$. 
To any $p\geq q$, we associate the subset $X^p_\bullet\subset [p+q]$ of cardinal $r=\sharp X_\bullet$ defined as 
$X^p_\bullet=X_+\coprod \{p+q+1-x,x\in X_-\}$.

The complement of $X^p_\bullet$, $(X^p_\bullet)^c=\{u_1<\cdots<u_{p+q-r}\}$, contains $(X^p_\bullet)^o=\{p+q+1-x,x\in X^p_\bullet\}$. Hence we can define 
the following subsets of cardinal $r$ : $(X^p_\bullet)'=\{i\in [p+q-r], u_i\in (X^p_\bullet)^o\}$ and 
$$
\widetilde{X}^p_\bullet=\{k\in [p+q-r], p+q-r+1-k\in (X^p_\bullet)'\}.
$$

Finally we associate the cycles classes $[\Xgot_{X^p_\bullet}]\in H^*(\G(r,p+q),\Z)$ and $[\Xgot_{\widetilde{X}^p_\bullet}]\in H^*(\G(r,p+q-r),\Z)$ to any  polarized subset $X_\bullet\subset [q]$ of cardinal $r\geq 1$, and to any $p\geq q$.

Recall that $\lambda(X^p_\bullet)$ (resp. $\lambda(\widetilde{X}^p_\bullet)$) is the partition associated to the subset $X^p_\bullet\subset [p+q]$ (resp. 
$\widetilde{X}^p_\bullet\subset [p+q-r]$). When $I_\bullet, J_\bullet,K_\bullet$ are three polarized subsets of $[q]$ of same cardinal $r$, the relation $[\Xgot_{I^p_\bullet}]\cdot[\Xgot_{J^p_\bullet}]\cdot[\Xgot_{K^p_\bullet}]=\ell [pt]\quad \ell\geq 1$, in $H^*(\G(r,p+q),\Z)$ is equivalent to $(\lambda(I^p_\bullet),\lambda(J^p_\bullet),\lambda(K^p_\bullet)-(p+q-r) 1^r)\in\horn(r)$. In the same way, the relation 
$[\Xgot_{\widetilde{I}^p_\bullet}]\cdot[\Xgot_{\widetilde{J}^p_\bullet}]\cdot[\Xgot_{\widetilde{K}^p_\bullet}]=\ell [pt]\quad \ell\geq 1$, in $H^*(\G(r,p+q-r),\Z)$ 
is equivalent to $(\lambda(\widetilde{I}^p_\bullet),\lambda(\widetilde{J}^p_\bullet),\lambda(\widetilde{K}^p_\bullet)-(p+q-2r) 1^r)\in\horn(r)$.

Let us rewrite a weak version of Theorem \ref{theo:singular-p-q-version-2} with these new notations.

\begin{theorem}\label{th:singular-horn}
An element $(x,y,z)\in(\R^q_{++})^3$ belongs to $\sing(p,q)$ if and only if 
$(\star)_{I_\pm,J_\pm,K_\pm}$ holds for any triplet $I_\bullet, J_\bullet,K_\bullet$ of polarized subsets of $[q]$ satisfying the following conditions~:
\begin{itemize}
\item[1.] $\sharp I_\bullet=\sharp J_\bullet=\sharp K_\bullet=r\in [q]$,
\item[2.] $(\lambda(I^p_\bullet),\lambda(J^p_\bullet),\lambda(K^p_\bullet)-(p+q-r) 1^r)\in\horn(r)$,
\item[3.] $(\lambda(\widetilde{I}^p_\bullet),\lambda(\widetilde{J}^p_\bullet),\lambda(\widetilde{K}^p_\bullet)-(p+q-2r) 1^r)\in\horn(r)$.
\end{itemize}
\end{theorem}

\medskip

Conditions {\em 2.} and {\em 3.} implies respectively that $|\lambda(I^p_\bullet)|+|\lambda(J^p_\bullet)|+|\lambda(K^p_\bullet)|=(p+q-r)r$ and 
$|\lambda(\widetilde{I}^p_\bullet)| + |\lambda(\widetilde{J}^p_\bullet)| + |\lambda(\widetilde{K}^p_\bullet)| = (p+q-2r)r$. We end this section by 
calculating precisely these two relations.

\medskip

\begin{definition}
For a polarized subset $X_\bullet\subset [q]$, we denote by $\delta_{X_\bullet}$ the cardinal of the set 
$\{(a,b)\in X_+\times X_-, a<b\}$.

Let $r\in [q]$. For any $0\leq \alpha \leq r$, we denote by $1_\alpha^r=(1,\ldots,1,0,\ldots,0)\in \mathbb{R}^r$ the vector where $1$ appears $\alpha$ times and $0$ appears 
$r-\alpha$ times.
\end{definition}

\begin{lemma}\label{lem:lambda-X}Let $X_\bullet\subset [q]$ be a polarized subset of cardinal $r\geq 1$. 
\begin{enumerate}
\item For any $p\geq q$, we have 
\begin{equation}\label{eq:codim-X}
|\lambda(X^p_\bullet)|= |X_-| - |X_+|+ (p+q+1)\sharp X_+ -\tfrac{r(r+1)}{2},
\end{equation}
\begin{equation}\label{eq:codim-X-tilde}
|\lambda(\widetilde{X}^p_\bullet)|= |X_-| - |X_+|+ (p+q+1)\sharp X_+ -\tfrac{r(r+1)}{2} - (\sharp X_+)^2 -2\delta_{X_\bullet}\cdot
\end{equation}
Hence $|\lambda(X^p_\bullet)|-|\lambda(\widetilde{X}^p_\bullet)|=(\sharp X_+)^2 +2\delta_{X_\bullet}$.
\item For any $p'\geq p\geq q$, we have
$$
\lambda(X^{p'}_\bullet)-\lambda(X^p_\bullet)= (p'-p) 1_{\sharp X_+}^r\quad {\rm and}\quad 
\lambda(\widetilde{X}^{p'}_\bullet)-\lambda(\widetilde{X}^p_\bullet)= (p'-p) 1_{\sharp X_+}^r\cdot
$$
\end{enumerate}
\end{lemma}

{\em Proof :} Let us prove the first point. A direct computation gives $|X^p_\bullet|=|X_+|-|X_-| + (p+q+1)(r-\sharp X_+)$ and then, as 
$|\lambda(X^p_\bullet)|=(p+q-r)r+\tfrac{r(r+1)}{2}-|X^p_\bullet|$, we obtain (\ref{eq:codim-X}).

Recall that $(X^p_\bullet)'\subset [p+q-r]$ is defined as the 
collection of the ranks of the elements of $\{p+q+1-x,x\in X^p_\bullet\}$ in $[n]-X^p_\bullet$. A short computation gives that 
\begin{equation}\label{eq:X-prime}
|(X^p_\bullet)'|= |X_-| - |X_+| -2\delta_{X_\bullet}+ (p+q+1-\sharp X_+)\sharp X_+.
\end{equation}
and $|\widetilde{X}^p_\bullet|=(p+q-r+1)r  - |(X^p_\bullet)'|$. 
As $|\lambda(\widetilde{X}^p_\bullet)|=(p+q-2r)r+\tfrac{r(r+1)}{2}-|\widetilde{X}^p_\bullet|$, we obtain
\begin{eqnarray*}
|\lambda(\widetilde{X}^p_\bullet)|&=&(p+q-2r)r +\tfrac{r(r+1)}{2}-(p+q-r+1)r +|(\widehat{X})'|\\
&=&|(X^p_\bullet)'|-\tfrac{r(r+1)}{2}\\
&=&|X_-| - |X_+|+ (p+q+1)\sharp X_+ -\tfrac{r(r+1)}{2} - (\sharp X_+)^2  -2\delta_{X_\bullet}\cdot
\end{eqnarray*}
The last equality is a consequence of (\ref{eq:X-prime}). 

The second point, which can be proved by direct calculation, is left to the reader. $\Box$

\medskip

\begin{corollary}\label{coro:inequality-p-q}
The inequality $(\star)_{I_\pm,J_\pm,K_\pm}$ appears in the description of the convex cone $\sing(p,q)$ only if 
the triplet $I_\bullet,J_\bullet,K_\bullet$ of polarized subsets of $[q]$ satisfies the relations : 
\begin{enumerate}
\item[(C1)] $\sharp I_\bullet =\sharp J_\bullet = \sharp K_\bullet =r\in [q]$,
\item[(C2)] $|I_+| + |J_+| + | K_+| - (|I_-| + |J_-| + K_-|)+\tfrac{r(r+1)}{2}=(p+q+1)(\sharp I_+ + \sharp J_+ + \sharp K_+ -r)$,
\item[(C3)] $(\sharp I_+)^2 + (\sharp J_+)^2 + (\sharp K_+)^2 + 2(\delta_{I_\bullet} +\delta_{J_\bullet} +\delta_{K_\bullet})= r^2$.
\end{enumerate} 
\end{corollary}

\medskip

\subsection{Examples : $\sing(p,1)$, $\sing(p,2)$ and $\sing(3,3)$}\label{sec:example-horn-1-2-3}

In the next examples, the notation $X_\bullet=\{i_1^\epsilon,\ldots, i_r^\epsilon\}$ means that $X_+=\{i_k^\epsilon, \epsilon=+\}$ and $X_-=\{i_k^\epsilon, \epsilon=-\}$.

\begin{example}\label{ex:singular-p-1}
Let $p\geq q=1$. Thus $(a,b,c)\in(\R_+)^3$ belongs to $\sing(p,1)$ if and only if the Weyl inequalities holds :
$a+b\geq c$, $a+c\geq b$, and $b+c\geq a$.
\end{example}

{\em Proof :} Let $I_\bullet,J_\bullet,K_\bullet\subset [1]$ satisfying the relations  of Theorem \ref{th:singular-horn}. If we use relation {\em (C3)} of Corollary 
\ref{coro:inequality-p-q}, we obtain that $(\sharp I_+)^2 + (\sharp J_+)^2 + (\sharp K_+)^2=1$. Up to a permutation, we must have $I_\bullet=J_\bullet=\{1^-\}$ and
$K_\bullet=\{1^+\}$. This case corresponds to the Weyl inequality $a+b\geq c$. $\Box$

\medskip

\begin{example}\label{ex:singular-p-2}
Let $p\geq q=2$. The element $(a,b,c)\in(\R^2_{++})^3$ belongs to $\sing(p,2)$ if and only if the following 18 inequalities holds 
\begin{enumerate}
\item the Weyl inequalities
\begin{itemize}
\item $a_1+ b_1\geq  c_1$ \quad (and 2 permutations),
\item $a_1+ b_2\geq  c_2$ \quad (and 5 permutations),
\end{itemize}
\item the Lidskii inequalities
\begin{itemize}
\item $a_1+a_2+b_1+b_2\geq c_1+c_2$\quad (and 2 permutations),
\end{itemize}
\item the signed Lidskii inequalities
\begin{itemize}
\item $a_1 + a_2 + b_1 - b_2 \geq c_1 - c_2$\quad (and 5 permutations).
\end{itemize}
\end{enumerate}
\end{example}

\medskip

{\em Proof :} Let $I_\bullet,J_\bullet,K_\bullet\subset [2]$ satisfying the relations  of Theorem \ref{th:singular-horn}. Let $r:=\sharp I_\bullet=\sharp J_\bullet=\sharp K_\bullet$.

First case : $r=1$.  Relation {\em (C3)} of Corollary \ref{coro:inequality-p-q} gives $(\sharp I_+)^2 + (\sharp J_+)^2 + (\sharp K_+)^2=1$. Up to a permutation, 
we have $I_{\bullet}=\{i^-\}$, $J_{\bullet}=\{j^-\}$, $K_{\bullet}=\{k^+\}$ with $i,j,k\in\{1,2\}$.   Relation {\em (C2)} of Corollary \ref{coro:inequality-p-q} 
imposes that $k=i+j-1$, thus the triplet $(i,j,k)$ belongs to $(1,1,1)$, $(1,2,2)$ and $(2,1,2)$. All this cases corresponds to the Weyl inequalities : $a_1+ b_1\geq  c_1$, $a_1+ b_2\geq  c_2$, and $a_2+ b_1\geq  c_2$.

Second case : $r=2$. Here relation {\em (C3)} becomes  
$$
(\sharp I_+)^2 + (\sharp J_+)^2 + (\sharp K_+)^2 + 2(\delta_{I_\pm}+\delta_{J_\pm}+\delta_{K_\pm})=4.
$$
Up to a permutation, we can suppose that $\sharp I_+\leq \sharp J_+\leq\sharp K_+$. There are two possibilities :
\begin{itemize}
\item $(\sharp I_+, \sharp J_+,\sharp K_+)=(0,0,2)$. Thus $I_\bullet=J_\bullet=\{1^-,2^-\}$ and $K_\bullet=\{1^+,2^+\}$. This case corresponds to the Lidskii inequality
$a_1+a_2+b_1+b_2\geq c_1+c_2$.
\item $(\sharp I_+, \sharp J_+,\sharp K_+)=(0,1,1)$. Up to a permutation, we have $\delta_{K_\bullet}=1$ and $\delta_{J_\bullet}=0$. Here $I_\bullet=\{1^-,2^-\}$,  
$J_\bullet=\{1^-,2^+\}$ and $K_\bullet=\{1^+,2^-\}$. This situation corresponds to the signed Lidskii inequality $a_1+a_2+b_1-b_2\geq c_1-c_2$. $\Box$
\end{itemize}

\medskip

\begin{example}\label{ex:singular-3-3}
An element $(a,b,c)\in(\R^3_{++})^3$ belongs to $\sing(3,3)$ if and only if the following 87 inequalities holds :
\begin{enumerate}
\item the Weyl inequalities
\begin{itemize}
\item $a_1+ b_1\geq  c_1$ \quad (and 2 permutations),
\item $a_1+ b_2\geq  c_2$ \quad (and 5 permutations),
\item $a_2 + b_2\geq c_3$ \quad (and 2 permutations),
\item $a_1 + b_3\geq c_3$ \quad (and 5 permutations).
\end{itemize}
\item the Lidskii inequalities
\begin{itemize}
\item $a_1 + a_2 + b_1 + b_2\geq c_1+c_2$ \quad (and 2 permutations),
\item $a_1 + a_2 + b_1 + b_3\geq c_1+c_3$ \quad (and 5 permutations),
\item $a_1 + a_2 + b_2 + b_3\geq c_2+c_3$ \quad (and 5 permutations),
\item $a_1 + a_2 + a_3 + b_1 + b_2 + b_3 \geq c_1 + c_2 + c_3$ \quad (and 2 permutations).
\end{itemize}
\item the signed Lidskii inequalities
\begin{itemize}
\item $a_1 + a_2 + b_1 - b_2 \geq c_1 - c_2$ \quad (and 5 permutations),
\item $a_1 + a_2 + b_1 - b_3 \geq c_1 - c_3$ \quad (and 5 permutations),
\item $a_1 + a_2 + b_2 - b_3 \geq c_2 - c_3$ \quad (and 5 permutations),
\item $a_1 + a_2 + a_3 + b_1 + b_2 - b_3 \geq c_1 + c_2 - c_3$ \quad (and 5 permutations),
\item $a_1 + a_2 + a_3 + b_1 - b_2 + b_3 \geq c_1 - c_2 + c_3$ \quad (and 5 permutations),
\item $a_1 + a_2 + a_3 - b_1 + b_2 + b_3 \geq - c_1 + c_2 + c_3$ \quad (and 5 permutations).
\end{itemize}
\item others inequalities
\begin{itemize}
\item $a_1 + a_3 + b_1 + b_3\geq c_2+c_3$ \quad (and 2 permutations),
\item $a_1 + a_3 + b_1 - b_3\geq c_2-c_3$ \quad (and 5 permutations),
\item $(a_1 + a_2 - a_3) + (b_1 - b_2 + b_3) + (-c_1 + c_2 + c_3)\geq 0$ \quad (and 5 permutations).
\end{itemize}
\end{enumerate}
\end{example}

{\em Proof :} Let $I_\bullet,J_\bullet,K_\bullet\subset [3]$ satisfying the relations  of Theorem \ref{th:singular-horn}. Let $r:=\sharp I_\bullet=\sharp J_\bullet=\sharp K_\bullet$.

\underline{First case : $r=1$.} 

We obtain here the Weyl inequalities $a_i+b_j\geq c_{i+j-1}$ (up to permutations).

\underline{Second case : $r=2$.} 

Let us use the relation 
$(\sharp I_+)^2 + (\sharp J_+)^2 + (\sharp K_+)^2 + 2(\delta_{I_\pm}+\delta_{J_\pm}+\delta_{K_\pm})=4$.
Up to a permutation, we can suppose that $\sharp I_+\leq \sharp J_+\leq\sharp K_+$. There are two possibilities :
\begin{itemize}
\item $(\sharp I_+, \sharp J_+,\sharp K_+)=(0,0,2)$. Here we have $I_+=J_+=K_-=\emptyset$. Relation {\em (C3)} of Corollary \ref{coro:inequality-p-q} gives 
$|I_-|+|J_-|=|K_+|+3$. Knowing that $|I_-|,|J_-|, |K_+|\in\{3,4,5\}$, we obtain few cases related to Lidskii inequalities:
\begin{enumerate}
\item $(|I_-|,|J_-|, |K_+|)=(3,3,3)$ :  $a_1 + a_2 + b_1 + b_2\geq c_1+c_2$,
\item $(|I_-|,|J_-|, |K_+|)=(3,4,4)\ or\  (4,3,4)$ : $a_1 + a_2 + b_1 + b_3\geq c_1+c_3$ or $a_1 + a_3 + b_1 + b_2\geq c_1+c_3$,
\item $(|I_-|,|J_-|, |K_+|)=(3,5,5)\ or \ (5,3,5)$ : $a_1 + a_2 + b_2 + b_3\geq c_2+c_3$ or $a_2 + a_3 + b_1 + b_2\geq c_2+c_3$.
\end{enumerate}
The only case left is $(|I_-|,|J_-|, |K_+|)=(4,4,5)$. Here $I_\bullet=J_\bullet=\{1^-,3^-\}$ and $K_\bullet=\{2^+,3^+\}$. A small computation shows 
that the Horn conditions are satisfied in this case. We obtain the inequality $a_1 + a_3 + b_1 + b_3\geq c_2+c_3$.

\item $(\sharp I_+, \sharp J_+,\sharp K_+)=(0,1,1)$. Up to a permutation, we can suppose that $\delta_{K_\bullet}=1$ and $\delta_{J_\bullet}=0$. Thus 
$I_\bullet=\{i_1^-<i_2^-\}$, $I_\bullet=\{j_1^-<j_2^+\}$ and $K_\bullet=\{k_1^+<k_2^-\}$. Relation {\em (C2)} of Corollary \ref{coro:inequality-p-q} gives 
$\underbrace{i_1^- + i_2^-}_\alpha + \underbrace{j_1^- + k_1^-}_\beta  =\underbrace{j_2^+ + k_1^+}_\gamma  + 3$.

Knowing that $\alpha,\beta,\gamma\in\{3,4,5\}$, we obtain 
\begin{enumerate}
\item $(\alpha,\beta,\gamma)= (3,3,3)$, $(3,4,4)$, or $(3,5,5)$. It gives us the signed Lidskii inequalities.
\item $(\alpha,\beta,\gamma)= (4,3,4)$, or $(5,3,5)$. No solutions.
\item $(\alpha,\beta,\gamma)= (4,4,5)$. Here $I_\bullet=\{1^-, 3^-\}$, $I_\bullet=\{1^-, 3^+\}$ and $K_\bullet=\{2^+,3^-\}$. A small computation 
shows that the Horn conditions are satisfied in this case. It corresponds to the inequality $a_1 + a_3 + b_1 - b_3\geq c_2 - c_3$.

\end{enumerate}
\end{itemize}

\medskip

\underline{Last case : $r=3$.} 

Let start with the relation {\em (C3)} :
$(\sharp I_+)^2 + (\sharp J_+)^2 + (\sharp K_+)^2 + 2(\delta_{I_\pm}+\delta_{J_\pm}+\delta_{K_\pm})=9$.
Up to a permutation, we can suppose that $\sharp I_+\leq \sharp J_+\leq\sharp K_+$. First, let us list the triples $(\sharp I_+\leq \sharp J_+\leq\sharp K_+)$ 
that cannot satisfy {\em (C3)} : $(0,0,0)$; $(0,0,1)$; $(0,1,1)$; $(0,0,2)$; $(1,1,2)$; $(0,2,2)$; $(2,2,2)$. There are four cases left :
\begin{enumerate}
\item $(\sharp I_+,\sharp J_+ ,\sharp K_+)=(0,0,3)$. Here $I_\bullet=J_\bullet= \{1^-, 2^-, 3^-\}$ and $K_\bullet= \{1^+, 2^+, 3^+\}$. We obtain the  Lidskii inequality 
$a_1 + a_2 + a_3 + b_1 + b_2 + b_3\geq c_1+ c_2 + c_3$.
\item $(\sharp I_+,\sharp J_+ ,\sharp K_+)=(0,1,2)$. We must have $\delta_{J_\bullet}+\delta_{K_\bullet}=2$, so we have three possibilities
\begin{itemize}
\item[a)] $(\delta_{J_\bullet},\delta_{K_\bullet})=(0,2)$ : here $I_\bullet= \{1^-, 2^-, 3^-\}$, $J_\bullet= \{1^-, 2^-, 3^+\}$ and $K_\bullet= \{1^+, 2^+, 3^-\}$. We obtain the signed 
Lidskii inequality $a_1 + a_2 + a_3 + b_1 + b_2 - b_3\geq c_1+ c_2 - c_3$.
\item[b)] $(\delta_{J_\bullet},\delta_{K_\bullet})=(1,1)$ : here $I_\bullet= \{1^-, 2^-, 3^-\}$, $J_\bullet= \{1^-, 2^+, 3^-\}$ and $K_\bullet= \{1^+, 2^-, 3^+\}$. We obtain the signed 
Lidskii inequality $a_1 + a_2 + a_3 + b_1 - b_2 + b_3\geq c_1- c_2 + c_3$.
\item[c)] $(\delta_{J_\bullet},\delta_{K_\bullet})=(2,0)$ : here $I_\bullet= \{1^-, 2^-, 3^-\}$, $J_\bullet= \{1^+, 2^-, 3^-\}$ and $K_\bullet= \{1^-, 2^+, 3^+\}$. We obtain the signed 
Lidskii inequality $a_1 + a_2 + a_3 - b_1 + b_2 + b_3\geq - c_1+ c_2 + c_3$.
\end{itemize}
\item $(\sharp I_+,\sharp J_+ ,\sharp K_+)=(1,2,2)$. We must have $\delta_{I_\bullet}=\delta_{J_\bullet}=\delta_{K_\bullet}=0$, thus $I_\bullet= \{1^-, 2^-, 3^+\}$, $J_\bullet= \{1^-, 2^+, 3^+\}$ and $K_\bullet= \{1^-, 2^+, 3^+\}$. The associated partitions are $\lambda(I^3_\bullet)=(1,0,0)$ and $\lambda(J^3_\bullet)=\lambda(K^3_\bullet)=(2,2,0)$. As $((1,0,0),(2,2,0),(-1,-1,-3))$ does not belongs to $\horn(3)$, this triplet 
$(I_\bullet,J_\bullet,K_\bullet)$ does not satisfies the conditions of Theorem \ref{th:singular-horn}.
\item $(\sharp I_+,\sharp J_+ ,\sharp K_+)=(1,1,1)$. We must have $\delta_{I_\bullet}+\delta_{J_\bullet}+\delta_{K_\bullet}=3$. Up to a permutation, two situations arise 
\begin{itemize}
\item[a)] $\delta_{I_\bullet}=\delta_{J_\bullet}=\delta_{K_\bullet}=1$: here $I_\bullet=J_\bullet= K_\bullet=\{1^-, 2^+, 3^-\}$. The corresponding partition is $\mu=(2,1,0)$, and one checks that 
$(\mu,\mu,\mu-3\cdot 1^3)\in \horn(3)$. Thus this triplet satisfies the conditions of Theorem \ref{th:singular-horn}. The corresponding inequality, 
$$
a_1+a_3+b_1+b_3+c_1+c_3\geq a_2+b_2+c_2,
$$
is nevertheless redundant since it is a consequence of the inequality $x_1\geq x_2$ satisfied by $a,b,c$.

\item[b)] $\delta_{I_\bullet}=0$, $\delta_{J_\bullet}=1$ and $\delta_{K_\bullet}=2$ :  here $I_\bullet= \{1^-, 2^-, 3^+\}$, $J_\bullet= \{1^-, 2^+, 3^-\}$ and $K_\bullet= \{1^+, 2^-, 3^-\}$. The corresponding partition are 
$\lambda=\lambda(I^3_\bullet)=(1,0,0)$, $\mu=\lambda(J^3_\bullet)=(2,1,0)$ and $\nu=\lambda(K^3_\bullet)=(3,1,1)$. One checks  that 
$(\lambda,\mu,\nu-3\cdot 1^3)\in \horn(3)$, thus the triplet $(I_\bullet,J_\bullet,K_\bullet)$ satisfies the conditions of Theorem \ref{th:singular-horn}. The corresponding inequality is $a_1+a_2+b_1+b_3+c_2+c_3\geq a_3+b_2+c_1$.

\end{itemize}

\end{enumerate}

\begin{remark}
The case where $I_\bullet,J_\bullet, K_\bullet$ are all equal to $A_\bullet=\{1^-, 2^+, 3^-\}$ yields the example of a triplet satisfying the Horn conditions but giving a 
redundant inequality. This redundancy phenomenon can be explained as the triplet $(A_\bullet ,A_\bullet,A_\bullet)$ does not satisfy the cohomological conditions of 
Theorem \ref{theo:singular-p-q-version-2} : one has $[\Xgot_{A^3_\bullet}]\cdot [\Xgot_{A^3_\bullet}]\cdot [\Xgot_{A^3_\bullet}]=2 [pt]$  in 
$H^*(\G(3,6),\Z)$ (see Example 4.2.1 in \cite{Berline-Vergne-Walter18} ).
\end{remark}

\subsection{The convex cone $\sing(\infty,q)$}\label{sec:infty-singular-horn}

Let $p\geq q$. To any matrix $X\in M_{p,q}(\C)$, we associate 
$$
X^\vee=\begin{pmatrix}X\\0\cdots 0\end{pmatrix}\in M_{p+1,q}(\C)
$$
It is immediate to see that $X$ and $X^\vee$ have the same singular values : $\tau(X)=\tau(X^\vee)$ for any $X\in M_{p,q}(\C)$.
From this basic consideration, we see that if $(\tau(A),\tau(B),\tau(A+B))\in \sing(p,q)$ then 
$(\tau(A),\tau(B),\tau(A+B))=(\tau(A^\vee),\tau(B^\vee),\tau(A^\vee + B^\vee))$ belongs to $\sing(p+1,q)$. Hence, 
$\sing(p,q)\subset  \sing(p+1,q)$ for any $p\geq q$. 

This leads us to define the following convex cone 
$$
\sing(\infty,q)=\bigcup_{p\geq q}\sing(p,q).
$$

Let us denote by $\Ecal_q$ the set formed by all inequalities 
$$
(\star)_{I_\pm,J_\pm,K_\pm} :\quad  \qquad|\,x\,|_{I_+}+|\,y\,|_{J_+} +|\, z\,|_{K_+}\leq |\,x\,|_{I_{-}}+|\,y\,|_{J_{-}} +|\, z\,|_{K_{-}},\qquad x,y,z\in \mathbb{R}^q.
$$
associated to triplet $(I_\bullet,J_\bullet,K_\bullet)$ of polarized subsets of $[q]$ of same cardinal.

\begin{definition}
Let us denote by $\Ecal_{p,q}\subset\Ecal_q$ the set of inequalities associated to triplet $(I_\bullet,J_\bullet,K_\bullet)$ satisfying the Horn conditions of 
Theorem \ref{th:singular-horn}. Let $\Ecal_{p,q}^{\rm min}\subset \Ecal_{p,q}$ be the minimal system of inequalities describing the convex cone $\sing(p,q)$.

An inequality $(\star)_{I_\pm,J_\pm,K_\pm}$ is called {\em regular} if there are twice as many terms on the right side as on the left side, e.g. $\sharp I_- +  \sharp J_- + \sharp K_- = 2(\sharp I_+ + \sharp J_+ + \sharp K_+)$.
\end{definition}

We start with a basic observation.
\begin{lemma}
There exists  $p_o\geq q$, such that $\sing(\infty,q)=\sing(p_o,q)$. In other words, 
$\sing(p,q)=\sing(p+1,q)$ for any $p\geq p_o$. 
\end{lemma}
{\em Proof :} Thanks to Theorem \ref{th:singular-horn}, we know that for any $p\geq q$ the convex cone $\sing(p,q)$ is determined by a finite system of 
inequalities $\Ecal_{p,q}\subset\Ecal_q$. As $\Ecal_q$ is finite, there exits $\Rcal\subset \Ecal_q$ such that $\N(\Rcal):=\{p\geq q, \Ecal_{p,q}=\Rcal\}$ is infinite. 
Take $p_o=\inf \N(\Rcal)$. Let us verify that $\sing(\infty,q)\subset\sing(p_o,q)$ : as the other inclusion $\sing(p_o,q)\subset\sing(\infty,q)$ is obvious, the proof 
will be completed. Let $x\in \sing(\infty,q)$ : there exists $p\geq q$ such that $x\in\sing(p,q)$. As $\N(\Rcal)$ is infinite, there exists $p'\in \N(\Rcal)$ such that 
$p'\geq p$. It follows that $x\in\sing(p,q)\subset \sing(p',q)=\sing(p_o,q)$. $\Box$

\medskip

We have now a nice characterization.

\begin{proposition}\label{prop:inequality-infty-q}
The identity $\sing(\infty,q)=\sing(p,q)$ holds if and only if all the inequalities of $\Ecal_{p,q}^{\rm min}$ are regular.
\end{proposition}

{\em Proof :} Let $p\geq q$ such that $\sing(\infty,q)=\sing(p,q)$. It means that $\Ecal_{p',q}^{\rm min}=\Ecal_{p,q}^{\rm min}$ for any $p'>p$. 
Let $(I_\bullet,J_\bullet,K_\bullet)$ such that $(\star)_{I_\pm,J_\pm,K_\pm}\in \Ecal_{p,q}^{\rm min}$. Corollary \ref{coro:inequality-p-q} tells us that 
$$
(C2)\qquad |I_+| + |J_+| + | K_+| - (|I_-| + |J_-| + K_-|)+\tfrac{r(r+1)}{2}=(p+q+1)(\sharp I_+ + \sharp J_+ + \sharp K_+ -r).
$$
Notice that $(\star)_{I_\pm,J_\pm,K_\pm}$ is a {\em regular} if and only if $\sharp I_+ + \sharp J_+ + \sharp K_+ -r=0$. Hence, if 
$(\star)_{I_\pm,J_\pm,K_\pm}$ is not regular, the identity {\em (C2)} does not hold anymore when we change $p$ in $p'>p$. It implies that 
$(\star)_{I_\pm,J_\pm,K_\pm}\notin \Ecal_{p',q}$ for any $p'>p$. It is in contradiction with the fact that $(\star)_{I_\pm,J_\pm,K_\pm}\in\Ecal_{p',q}^{\rm min}$.

Suppose now that all the inequalities of $\Ecal_{p,q}^{\rm min}$ are regular. In order to prove that $\sing(\infty,q)$ coincides with $\sing(p,q)$, it is sufficient to show that 
$\Ecal_{p,q}^{\rm min}\subset \Ecal_{p',q}$ for any $p'\geq p$. Let $(I_\bullet,J_\bullet,K_\bullet)$ such that $(\star)_{I_\pm,J_\pm,K_\pm}\in \Ecal_{p,q}^{\rm min}$. By definition, 
we have 
\begin{itemize}
\item[1.] $\sharp I_\bullet=\sharp J_\bullet=\sharp K_\bullet=r\in [q]$,
\item[2.] $\sharp I_+ +\sharp J_+ +\sharp K_+=r$,
\item[3.] $\Lambda^{p}:=(\lambda(I^p_\bullet),\lambda(J^p_\bullet),\lambda(K^p_\bullet)-(p+q-r) 1^r)\in\horn(r)$,
\item[4.] $\widetilde{\Lambda}^p:=(\lambda(\widetilde{I}^p_\bullet),\lambda(\widetilde{J}^p_\bullet),\lambda(\widetilde{K}^p_\bullet)-(p+q-2r) 1^r)\in\horn(r)$.
\end{itemize}
Identity {\em 2.} is due to our assumption that $(\star)_{I_\pm,J_\pm,K_\pm}$ is a regular inequality.

Our goal is achieved if we show that $\Lambda^{p'}$ and $\widetilde{\Lambda}^{p'}$ belong to  $\horn(r)$ for any $p'>p$.
If $0\leq \alpha\leq r$, we denote by $1_\alpha^r=(1,\ldots,1,0,\ldots,0)\in \mathbb{R}^r$ the element where $1$ appears $\alpha$-times. 
We have seen in Lemma \ref{lem:lambda-X} that for any polarized subset $X_\bullet\subset [q]$ of cardinal $r\geq 1$ we have 
$\lambda(X^{p'}_\bullet)-\lambda(X^p_\bullet)= \lambda(\widetilde{X}^{p'}_\bullet)-\lambda(\widetilde{X}^p_\bullet)= (p'-p) 1^r_{\sharp X_+}$, $\forall p'\geq p$. 
This leads to the relations 
$$
\Lambda^{p'}-\Lambda^{p}=\widetilde{\Lambda}^{p'}-\widetilde{\Lambda}^p=(p'-p)\left(1^r_{\sharp I_+},1^r_{\sharp J_+},1^r_{\sharp K_+}-1^{r}\right),\qquad \forall p'\geq p.
$$
Now it is an easy matter to check that $(1^r_{\sharp I_+},1^r_{\sharp J_+},1^r_{\sharp K_+}-1^{r})\in \horn(r)$ since $\sharp I_+ +\sharp J_+ +\sharp K_+=r$. 
Using the fact that $\horn(r)$ is a convex cone, we can conclude that $\Lambda^{p'}=\Lambda^p+(p'-p)(1^r_{\sharp I_+},1^r_{\sharp J_+},1^r_{\sharp K_+}-1^{r})$ and 
$\widetilde{\Lambda}^{p'}=\widetilde{\Lambda}^p+(p'-p)(1^r_{\sharp I_+},1^r_{\sharp J_+},1^r_{\sharp K_+}-1^{r})$ belongs to $\horn(r)$ for any $p'\geq p$. $\Box$

\medskip

\begin{corollary}
We have the following identities :
\begin{itemize}
\item  $\sing(\infty,1)=\sing(1,1)$,
\item $\sing(\infty,2)=\sing(2,2)$,
\item $\sing(\infty,3)=\sing(3,3)$.
\end{itemize}
\end{corollary}

{\em Proof :} The first two relations follow from the computations done in Examples \ref{ex:singular-p-1} and \ref{ex:singular-p-2}. In Example \ref{ex:singular-3-3}, we have proved that $\sing(3,3)$ is determined by a system of {\em regular} inequalities. The relation $\sing(\infty,3)=\sing(3,3)$ is then a consequence of the previous proposition. $\Box$

\subsection{Proof of Theorem \ref{theo:singular-p-q}}\label{sec:preuve-singular-horn}

We work with the conjugate linear involution $\sigma_{p,q}(g)=I_{p,q}(g^*)^{-1}I_{p,q}$ on $GL_n(\C)$, 
where $I_{p,q}$ is the diagonal matrix ${\rm Diag}(I_p,-I_q)$. Let $U_\C=GL_n(\C)$ embedded diagonally in $\tU_\C=GL_n(\C)\times GL_n(\C)$. 
We define on $\tU_\C$ the conjugate linear involution $\sigma_{p,q}(g_1,g_2)=(\sigma_{p,q}(g_1),\sigma_{p,q}(g_2))$.

We see then that $G=U(p,q)$ and $\tG=G\times G$ are the subgroups fixed by $\sigma_{p,q}$. The Cartan decomposition of the Lie algebra of 
$G$ is $\ggot=\kgot\oplus\pgot$ where $\kgot=\ugot(p)\times\ugot(q)$ and $\pgot\simeq M_{p,q}(\C)$ as a $U(p)\times U(q)$-module. 

Let $\tgot\subset \ugot(n)$ be the maximal torus composed by all the matrices of the form 
$$
i\,\begin{pmatrix} D(a)&  0& M(c)\\
0& D(b)& 0\\
 M(\tilde{c})&  0& D(\tilde{a})
\end{pmatrix},\qquad (a,b,c)\in \R^q\times \R^{p-q}\times \R^q.
$$
Here 
$M(c)=\begin{pmatrix}
0&\cdots & c_1\\
\vdots &\reflectbox{$\ddots$}& \vdots\\
c_q&\cdots & 0\\
\end{pmatrix}$ and $D(a)=\begin{pmatrix}
a_1&\cdots & 0\\
\vdots &\ddots& \vdots\\
0&\cdots & a_q\\
\end{pmatrix}$. We denoted $\tilde{x}=(x_q,\ldots,x_1)$ for any $x=(x_1,\ldots,x_q)\in\R^q$.

We notice that $\tgot$ is invariant under the involution $\sigma_{p,q}$, and that $\agot=i \tgot^{-\sigma_{p,q}}$ is the maximal abelian subspace of $\pgot$ consisting of all matrices of the form 
\begin{equation}\label{eq:maxabelian}
X(c):=\begin{pmatrix}0&  0& M(c)\\
0&  0& 0\\
M(\tilde{c}) & 0 & 0
\end{pmatrix},\qquad c\in\R^q.
\end{equation}

Let $f_i$ be the member of $\agot^*$ whose value on the matrix $X(a)$ is $a_i$. 
The set $\Sigma(\ggot)$ of restricted roots includes all linear functionals $\pm f_i\pm f_j$ with $i\neq j$ and $\pm 2 f_i$ for all $i$. Also the $\pm f_i$ 
are restricted roots if $p\neq q$ (see Chapter 6 in \cite{Knapp}). Thus we can choose the following Weyl chamber 
$$
\agot_+:=\{X(a), a_1\geq a_2\geq \cdots\geq a_q\geq 0\}.
$$
In this case the restricted Weyl group $W_\agot$ is the semi-direct product $\Sgot_q\rtimes \{1,-1\}^q$ where $\Sgot_q$ acts on $\R^q$ by permuting the indices and 
$\{1,-1\}^q=\{\epsilon=(\epsilon_1,\cdots,\epsilon_q)\ \vert\  \epsilon_k=\pm 1\}$ acts by sign changes : $\epsilon \cdot a= (\epsilon_1 a_1,\ldots, \epsilon_q a_q)$.

We start by determining the admissible elements.

\begin{lemma}\label{lem:admissible-singular-1}
Let $G=U(p,q)$ and $\tG=G\times G$. The {\em admissible} elements relatively to the $\tG\times G$-action on $\tG$ are\footnote{Up to multiplication by a positive rational.} of the form 
$(w_1\gamma_r,w_2\gamma_r,\gamma_3\gamma_r)$ where $w_i\in W_\agot,\forall i$ and 
$$
\gamma_r:=\begin{pmatrix}0&  0& M(c_r)\\
0&  0& 0\\
M(\tilde{c_r}) & 0 & 0
\end{pmatrix},
$$
where $c_r=(\underbrace{-1,\ldots,-1}_{ r\ times},0,\ldots,0)$ for some $r\in\{1,\ldots,q\}$.
\end{lemma}

{\em Proof :} First we note that $\tggot/\ggot\simeq \ggot$ as a $\ggot$-module, thus $\Sigma(\tggot/\ggot)=\Sigma(\ggot)$. We look at a rational vector 
$\gamma$ such that $(\star)\ {\rm Vect}\big(\Sigma(\ggot)\cap \gamma^\perp\big)={\rm Vect}\big(\Sigma(\ggot)\big)\cap \gamma^\perp$. Up to the action of 
the Weyl group $W_\agot$, we can suppose that $\gamma=(\gamma_1\geq\cdots\geq\gamma_q\geq 0)$ is dominant. Since the set $\Sigma(\ggot)$ is generated 
by the simple roots ${\rm S}(\ggot)=\{f_{i}-f_{i+1}, 1\leq i\leq q-1\}\cup\{f_q\}$, we see that for $\gamma$ dominant, the condition $(\star)$ is equivalent to asking that all linear forms of ${\rm S}(\ggot)$ except one must vanish against $\gamma$. This last condition implies that $\gamma$ is proportional to a certain $\gamma_r$ for 
$r\in\{1,\ldots,q\}$. $\Box$

\medskip

Now we have to detailed condition $b)$ of Theorem \ref{theo:version-levi-movable} for an admissible element 
$(w_1\gamma_r,w_2\gamma_r,w_3\gamma_r)$. To do so, we use the conjugation by the element $\theta\in O(n)$ defined by 
$$
\theta=\frac{1}{\sqrt{2}}
\begin{pmatrix}
I_q&0& J_q\\
0&\sqrt{2}I_{p-q} & 0\\
J_q& 0 & -I_q
\end{pmatrix}.
$$
Here $J_q=\begin{pmatrix}0&\cdots & 1\\ 0&\reflectbox{$\ddots$}& 0\\1& \cdots & 0\end{pmatrix}$ is a $q\times q$ matrix. We see that $\theta^2=I_n$ and that 
$\tilde{\tgot}:=Ad(\theta)(\tgot)\subset \ugot(n)$ is the (usual) maximal torus of diagonal matrices. More precisely, the image of the 
matrix $
\begin{pmatrix} D(a)&  0& M(c)\\
0& D(b)& 0\\
 M(\tilde{c})&  0& D(\tilde{a})
\end{pmatrix}
$ 
by $Ad(\theta)$ is equal to 
$
\begin{pmatrix} D(a+c)&  0&0\\
0& D(b)& 0\\
0&  0& D(\tilde{a}-\tilde{c})
\end{pmatrix}.
$

Let $E_\theta=\theta(\C^r)\subset F_\theta=\theta(\C^{n-r})$ be the image of the flag $\C^r\subset \C^{n-r}$ by $\theta$. The map 
$\varphi([g])= g(E_\theta\subset F_\theta)$ defines a diffeomorphism between the flag manifold 
$\Fcal_{\gamma_r}:=GL_n(\C)/\Pbb_{\gamma_r}$ and the two-steps flag variety $\F(r,n-r;n)$. Moreover, for any $w\in W_\agot$ the image of the Schubert variety 
$\Xgot_{w,\gamma_r}$ through the diffeomorphism $\varphi$ is equal to some $\theta(\Xgot_{\widehat{X}})\subset \F(r,n-r;n)$ where $X_+\cup X_-=w(\{1,\cdots, r\})$. To be more precise the polarized subset $X_+\cup X_-$ is equal to 
$\{|w(i)|, 1\leq i\leq r\}$ and $k\in X_+$ if and only if there exists $i\in[r]$ with $k=w(i)=|w(i)|$.

Through the diffeomorphism $\varphi$, we see that the condition 
$$
b)\qquad [\Xgot_{w_1,\gamma_r}]\odot_0 [\Xgot_{w_2,\gamma_r}]\odot_0 [\Xgot_{w_3,\gamma_r}]= [pt]\quad {\rm in}\quad H^*(\Fcal_{\gamma_r},\Z)
$$
becomes $[\Xgot_{\widehat{I}}]\odot_0 [\Xgot_{\widehat{J}}]\odot_0 [\Xgot_{\widehat{K}}]=[pt]$  in $H^*(\F(r,n-r;n),\Z)$, where $I_+\cup I_-=w_1(\{1,\cdots, r\})$, 
$J_+\cup J_-=w_2(\{1,\cdots, r\})$ and $K_+\cup K_-=w_3(\{1,\cdots, r\})$ are polarized subsets of $[q]$.


\end{document}